# Fast and stable global interpolation based on equidistant points

Xu-Qing Liu[a], Hao Liu[b], Jian-Ying Rong[c]

[a]*Faculty of Mathematics and Physics, Huaiyin Institute of Technology, Huai'an, 223003, China*
[b]*Department of Electronic Engineering, Tsinghua University, Beijing, 100084, China*
[c]*Jiangsu Vocational College of Electronics and Information, Huai'an, 223003, China*

**Abstract**

This paper presents the *symmetric wave interpolation* method for *stable global interpolation* using readily available *equidistant points*. Its key achievement is the integration of the practical utility of such points with the numerical stability of Chebyshev interpolation. Experimental results demonstrate that symmetric wave interpolation effectively suppresses the Runge phenomenon and, crucially, delivers accuracy that matches or even surpasses Chebyshev interpolation. This work thereby provides a robust and practical solution that bridges the long-standing gap between point accessibility and numerical stability in global interpolation.

*Keywords:* Global interpolation, Symmetric wave interpolation, Equidistant points, Chebyshev interpolation, Chebyshev points, (Barycentric) Lagrange interpolation, Runge phenomenon

*2020 MSC:* 65D05

## 1. Introduction

Consider the problem of constructing a *global interpolant* from data sampled on $n+1$ *equidistant points*. The classical Lagrange interpolation polynomial is the immediate theoretical solution [1]. However, it suffers from three significant drawbacks for practical computation [2, 3]: Its computational complexity is of $O(n^2)$, making it inefficient for a large $n$; It is also numerically unstable, and lacks extensibility, as adding a new point requires restarting the entire computation from scratch.

The barycentric Lagrange interpolation formula was introduced as a remedy [2]. It offers a compelling alternative with $O(n)$ complexity after weight precomputation, superior numerical stability, and extensibility similar to the Newton interpolation. This formulation rightly raised the hope of a more robust and efficient framework for global interpolation.

Despite these computational advantages, a fundamental barrier still remains: the well-known *Runge phenomenon* [4]. This term describes the severe oscillatory behavior and drastic error increase near the two endpoints of the interval when performing a global high-degree polynomial interpolation on equidistant points. Figure 1 illustrates this by interpolating the Runge function, $f_1(x) \triangleq 1/(1+25x^2)$. This is regarded as an inherent mathematical issue of the underlying polynomial, not a numerical effect of a particular formula [5]. Therefore, both the classical and the barycentric forms are equally afflicted when the points are equidistant.

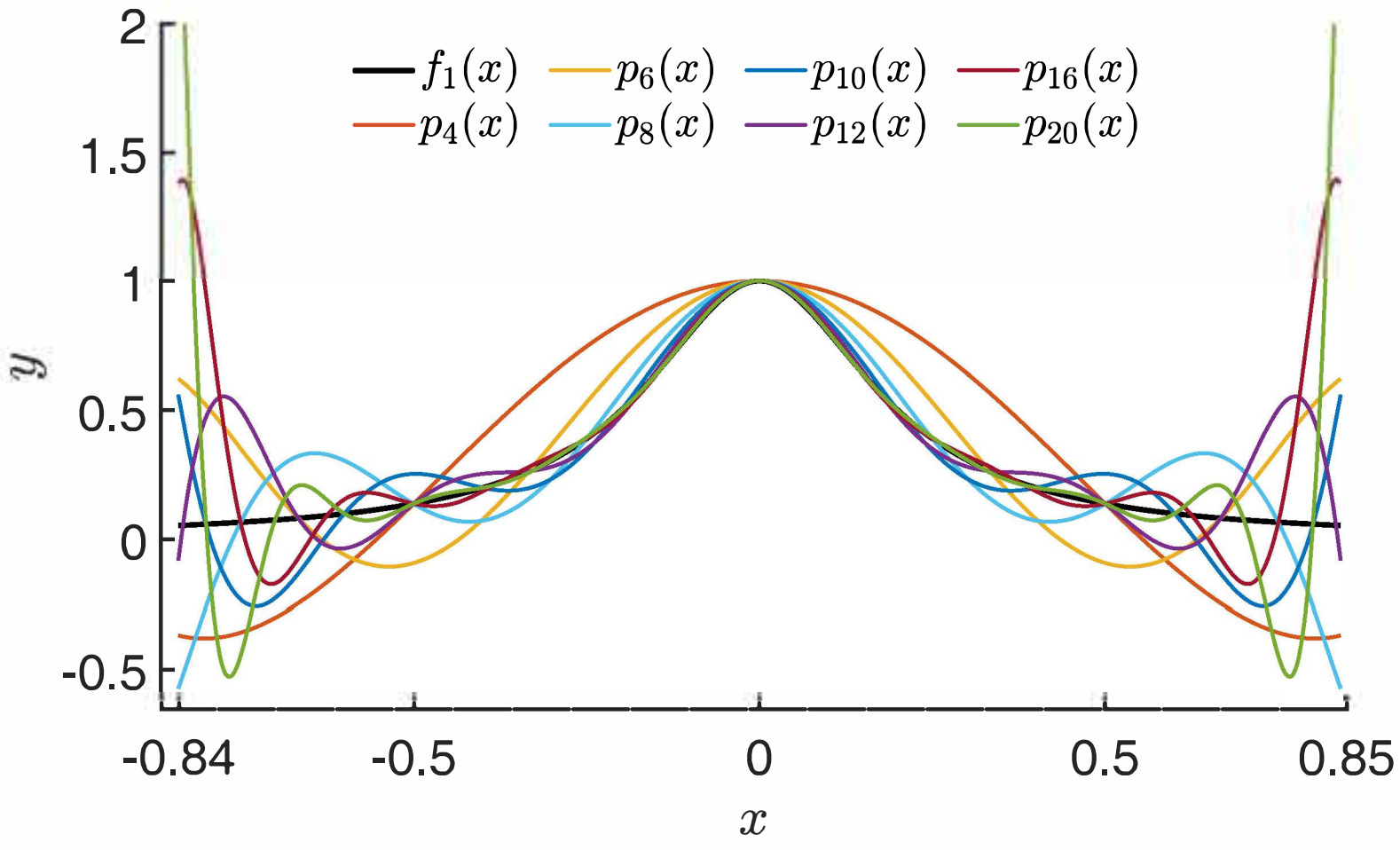

Figure 1: Polynomial interpolation of $f_1(x) \triangleq 1/(1+25x^2)$ on the $(n+1)$ equidistant points for $n = 4,\ 6,\ 8,\ 10,\ 12,\ 16,\ 20$.

To avoid the Runge phenomenon, existing research mainly follows two different paths, both with inherent limitations. The first path gives up the requirement for equidistant points. Instead, it uses a non-uniform point distribution such as the *Chebyshev points*, whose point-clustering property theoretically ensures a well-behaved and stable interpolation process, typically yielding flawless, exponentially convergent results [6, 7]. However, this approach requires free choice of point locations. This key condition cannot be met when working with the *most commonly used* fixed equidistant data. The second keeps equidistant points but abandons global interpolation. It includes piecewise low-degree polynomials like splines [8] or regularized fitting methods [2, 9–11]. The former loses global smoothness and the simplicity of analytical expression, while the latter violates the fundamental definition of interpolation because it does not pass exactly through all data points.

This reveals a long-standing challenge in numerical analysis: Does there exist a strict global interpolation method that works directly on given equidistant points? If this is the case, can it follow the exact interpolation condition and also effectively suppress the Runge phenomenon? This paper aims to address this challenge and fill the existing research gap. According to this conception, we will propose a new interpolation framework by transforming equidistant points to Chebyshev points via suitable mappings. The two-step framework, *transformation* followed by *interpolation*, is expected to inherit the superior qualities of Chebyshev interpolation and well work on the points that are equally spaced. In this sense, our method offers a powerful solution to this long-standing challenge in global interpolation with equidistant points.

The remainder of this paper is organized as follows: Section 2 provides a concise background on the relevant interpolation methods. Section 3 states the construction of our method, gives the truncation error bounds, and shows how it works. Section 4 presents a comprehensive numerical experiment, comparing the performance of our method based on equidistant points against the Chebyshev interpolation which, however, works only on Chebyshev points. Finally, Section 5 concludes the paper and outlines potential future research directions.

## 2. Background

This section begins with a concise review of the two formulations of Lagrange interpolation (classical and barycentric) and Chebyshev interpolation [12, 13]. It then explains why these methods typically exhibit numerical stability when applied to Chebyshev points, and how this computational advantage vanishes when using equidistant points.

### *2.1. Lagrange interpolation and barycentric Lagrange interpolation*

Let $f(x)$ be a function defined on the interval $[-1,1]$, sampled at $n+1$ equidistant points

$$x_i = -1 + \frac{2i}{n} \quad (i = 0, 1, \cdots, n).$$

The goal of global polynomial interpolation is to find an $n$-th degree polynomial $p_n(x)$ such that $p_n(x_i) = f(x_i) \triangleq y_i$ holds for all $i$. For the case of interpolating $\tilde{f}(\tilde{x})$ for $\tilde{x} \in [a,b]$, it coincides with interpolating $f(x)$ for $x \in [-1,1]$, by taking $x = [2\tilde{x} - (a+b)]/(b-a)$ and writing

$$\tilde{f}(\tilde{x}) = \tilde{f}\left(\frac{b-a}{2}x + \frac{a+b}{2}\right) \triangleq f(x) = f\left(\frac{2\tilde{x} - (a+b)}{b-a}\right).$$

In this sense, we only consider interpolating $f(x)$ for $x \in [-1,1]$, for convenience.

The classical Lagrange form can be expressed as

$$p_n(x) = \sum_{i=0}^{n} y_i \ell_i(x), \tag{2.1}$$

where $\ell_i(x) = \prod_{j=0,\, j\neq i}^{n} \frac{x-x_j}{x_i-x_j}$ $(i = 0, 1, \cdots, n)$ are the Lagrange basis polynomials. This form is powerful for theoretical analysis; however, it suffers from a computational complexity of $O(n^2)$, numerical instability, and poor extensibility.

The barycentric Lagrange interpolation reformulates the expression of (2.1) by introducing barycentric weights $w_i$, yielding the following barycentric formula:

$$p_n(x) = \frac{\sum_{i=0}^{n} \frac{w_i}{x-x_i} y_i}{\sum_{i=0}^{n} \frac{w_i}{x-x_i}}, \quad \text{with} \quad w_i = \frac{1}{\prod_{j=0,\, j\neq i}^{n} (x_i - x_j)}. \tag{2.2}$$

This form is mathematically equivalent to (2.1), while reducing the computational complexity to $O(n)$ after precomputing the weights. It also possesses superior numerical stability and is highly extensible.

### *2.2. Chebyshev interpolation*

Although barycentric Lagrange interpolation, $p_n(x)$ involved in (2.2), is numerically stable, its stability relies on a prerequisite: the points $x_0, x_1, \cdots, x_n$ must yield well-behaved weights $w_0, w_1, \cdots, w_n$. When the points are equidistantly distributed, the cancellation of large numbers drastically amplifies rounding errors, rendering the weights ill-conditioned and the interpolation result (especially near the interval boundaries) highly sensitive to computational noise, appearing as the Runge phenomenon. In contrast, Chebyshev points yield well-behaved barycentric weights, thereby theoretically guaranteeing numerical stability as well as uniform convergence [6].

Denote the *Chebyshev polynomials* and *Chebyshev points* (of the first and second kinds) by

$$
\begin{aligned}
T_k(z) &= \cos(k \arccos z), \ \ k = 0, 1, \cdots, n \\
z_i^{[1]} &= \cos\left(\frac{2i+1}{2(n+1)}\pi\right), \ \ i = 0, 1, \cdots, n \\
z_i^{[2]} &= \cos\left(\frac{i}{n}\pi\right), \ \ i = 0, 1, \cdots, n
\end{aligned}
$$

respectively. Then, the Chebyshev interpolation polynomials with respect to the first and second Chebyshev points can be written as

$$
p_n^{[1]}(z) = \sum_{k=0}^{n} \delta_k^{[1]} c_k^{[1]} T_k(z), \quad \text{with} \quad \delta_k^{[1]} = \begin{cases} 1/2, & k = 0 \\ 1, & k = 1, \cdots, n \end{cases} \tag{2.3}
$$

$$
p_n^{[2]}(z) = \sum_{k=0}^{n} \delta_k^{[2]} c_k^{[2]} T_k(z), \quad \text{with} \quad \delta_k^{[2]} = \begin{cases} 1/2, & k = 0, n \\ 1, & k = 1, \cdots, n-1 \end{cases} \tag{2.4}
$$

in which $c_k^{[1]}$ and $c_k^{[2]}$ are defined as follows:

$$
\begin{aligned}
c_k^{[1]} &= \frac{2}{n+1} \sum_{i=0}^{n} y_i T_k\left(z_i^{[1]}\right) = \frac{2}{n+1} \sum_{i=0}^{n} y_i \cos\left(\frac{k(2i+1)}{2(n+1)}\pi\right), \\
c_k^{[2]} &= \frac{2}{n} \sum_{i=0}^{n} \delta_i^{[2]} y_i T_k\left(z_i^{[2]}\right) = \frac{2}{n} \sum_{i=0}^{n} \delta_i^{[2]} y_i \cos\left(\frac{ki}{n}\pi\right).
\end{aligned}
$$

See [13, chatper 6] for the details. Note that the barycentric weights based on the Chebyshev points of the first and second kinds are of the closed forms [2]:

$$
w_i^{[1]} = (-1)^i \sin\left(\frac{(2i+1)\pi}{2(n+1)}\right), \quad w_i^{[2]} = (-1)^i \delta_i^{[2]}
$$

by cancelling all the factors independent of $i$. This indicates these weights are well-conditioned, ensuring the numerical stability and uniform convergence of interpolation.

Figure 2 presents the Chebyshev interpolation of $f_1(x)$ based on Chebyshev points of the first and second kinds, demonstrating that the Runge phenomenon can be significantly alleviated compared with the interpolation on equidistant points presented in Figure 1.

### *2.3. Motivation*

Although the interpolation using Chebyshev points is highly successful, it faces a practical limitation: *acquiring such specially distributed points is often much less convenient* than using equidistant points. For example, those data from many sources, such as sensor readings and historical financial records, are typically sampled at equidistant points, rather than at designed special points. Consequently, barycentric Lagrange interpolation or Chebyshev interpolation on Chebyshev points is high-performing but impractical for widespread use, especially for direct application in analyzing and modeling such data.

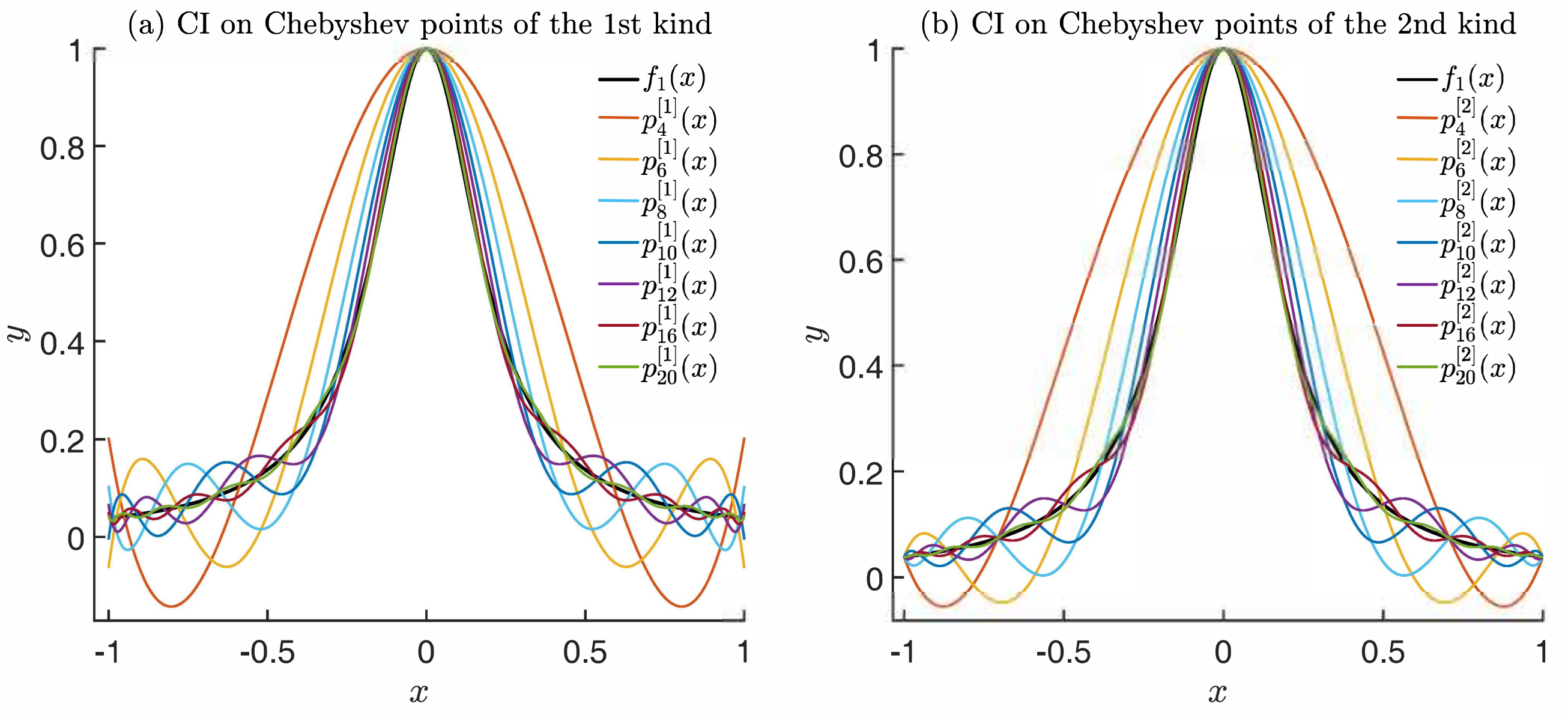

Figure 2: Chebyshev interpolation (CI) of $f_1(x)$ on the $(n+1)$ Chebyshev points of the first and second kinds for $n = 4, 6, 8, 10, 12, 16, 20$.

This motivates a fundamental question: can we adapt the superior qualities of Chebyshev interpolation to work effectively on equidistant points? The answer is "yes", noting that both equidistant points and Chebyshev points are regular and can be mapped to each other through appropriate transformations. Although this idea is still preliminary and there may be unknown challenges, we believe it is overall sound. In this sense, it allows a wide range of fields, even those without specialized sampling capabilities, to benefit from the accuracy and stability of *Chebyshev points based* Chebyshev interpolation.

## 3. Methodology

Motivated by the issue revealed in Section 2, we are led to the following insight: in order to combine the *stability of Chebyshev interpolation* (or the barycentric Lagrange interpolation) with the *convenience of equidistant points*, we can construct mappings that convert the distribution of the points from equidistant to Chebyshevian.

### *3.1. Mappings from equidistant to Chebyshevian*

Let $\tau_1$ and $\tau_2$ be two mappings defined as

$$\tau_1(z) \triangleq -\frac{2(n+1)}{n\pi} \arcsin z \quad \Rightarrow \quad \kappa_1(x) \triangleq \tau_1^{-1}(x) = -\sin\left(\frac{n\pi}{2(n+1)}x\right) \tag{3.1}$$

$$\tau_2(z) \triangleq -\frac{2}{\pi} \arcsin z \quad \Rightarrow \quad \kappa_2(x) \triangleq \tau_2^{-1}(x) = -\sin\left(\frac{\pi}{2}x\right) \tag{3.2}$$

respectively, where $\kappa_j$ is the inverse mapping of $\tau_j$ for $j = 1, 2$. Then, we have

$$x_i = -1 + \frac{2i}{n} \;\Rightarrow\; \begin{cases} \kappa_1(x_i) = \cos\left(\frac{2i+1}{2(n+1)}\pi\right) = z_i^{[1]} \;\Rightarrow\; x_i = \tau_1\left(z_i^{[1]}\right), \\ \kappa_2(x_i) = \cos\left(\frac{i}{n}\pi\right) \qquad\;\; = z_i^{[2]} \;\Rightarrow\; x_i = \tau_2\left(z_i^{[2]}\right). \end{cases}$$

That is, the mapping $\kappa_1$ defined by (3.1) converts the equidistant points into the Chebyshev points of the first kind, while the mapping $\kappa_2$ defined by (3.2) converts them into the second kind. These two mappings allow us to use Chebyshev points indirectly even if we only have equidistant points. The method is simple and easy to implement, it can make the interpolation converge fast (as the number of points increases) and remain numerically stable with respect to small data perturbations.

Here, it is mentioned that $\tau_2(z)$ can be viewed as a special form of the transformation of Kosloff and Tal-Ezer [14], except for a negative sign, by taking the parameter $\alpha \to 1$, if we only consider the form of the expression without regard to specific contexts.

*3.2. Interpolation formulae*

Inserting $x = \tau_j(z)$ into the function $f(x)$, we obtain $f(x) = f(\tau_j(z)) \triangleq g^{[j]}(z)$. This means *interpolating $f(x)$ on equidistant points* amounts to *interpolating $g^{[j]}(z)$ on Chebyshev points* of the $j$-th kind for $j = 1, 2$. Suppose the $n$-th derivative, namely $f^{(n)}(x)$, of $f(x)$ is continuous on $[-1, 1]$, and $f^{(n+1)}(x)$ exists on $(-1, 1)$. Then, we have the following theorems:

**Theorem 1.** *Let $p_n^{[1]}(z)$ be the Chebyshev interpolation polynomial, presented in (2.3), of $g^{[1]}(z)$ on $z_0^{[1]}, z_1^{[1]}, \cdots, z_n^{[1]}$. Then, the function interpolating $f(x)$ based on the equidistant points, namely $x_0, x_1, \cdots, x_n$, can be expressed as*

$$p_n^{[1]}(\kappa_1(x)) = p_n^{[1]}\left(-\sin\left(\frac{n\pi}{2(n+1)}x\right)\right) \triangleq q_n^{[1]}(x).$$

*Further, the truncation error bound is*

$$\max_{-1 \leqslant x \leqslant 1} \left|f(x) - q_n^{[1]}(x)\right| \leqslant \frac{2^{-n}}{(n+1)!} \max_{\kappa_1(1) < z < \kappa_1(-1)} \left|\left(g^{[1]}(z)\right)^{(n+1)}\right|.$$

*Proof.* Substituting $z = \kappa_1(x)$ into $p_n^{[1]}(z)$, we get $q_n^{[1]}(x)$ immediately. The truncation error bound can be derived similarly to the established error theory for Chebyshev interpolation in the literature, e.g., [15]. □

**Theorem 2.** *Let $p_n^{[2]}(z)$ be the Chebyshev interpolation polynomial, presented in (2.4), of $g^{[2]}(z)$ on $z_0^{[2]}, z_1^{[2]}, \cdots, z_n^{[2]}$. Then, the function interpolating $f(x)$ based on the equidistant points, namely $x_0, x_1, \cdots, x_n$, can be expressed as*

$$p_n^{[2]}(\kappa_2(x)) = p_n^{[2]}\left(-\sin\left(\frac{\pi}{2}x\right)\right) \triangleq q_n^{[2]}(x).$$

*Further, for any given $\varepsilon > 0$, exist $\lambda_1 \in (-1, 0)$ and $\lambda_2 \in (0, 1)$ such that*

$$\max_{-1 \leqslant x \leqslant 1} \left|f(x) - q_n^{[2]}(x)\right| \leqslant \max\left\{\varepsilon, \frac{2^{1-n}}{(n+1)!} \max_{\lambda_1 < z < \lambda_2} \left|\left(g^{[2]}(z)\right)^{(n+1)}\right|\right\}. \tag{3.3}$$

*Proof.* Before proving Theorem 2, note that the derivative of $\tau_2$ becomes infinite at $\pm 1$. This is why Theorem 2 handles the neighborhoods of the endpoints separately to estimate the remainder term. In addition, we note that the error bound for $q_n^{[2]}$ is twice as large as that for $q_n^{[1]}$.

Substituting $z = \kappa_2(x)$ into $p_n^{[2]}(z)$ gives $q_n^{[2]}(x)$. Note that $h(x) \triangleq f(x) - q_n^{[2]}(x)$ is continuous at $\pm 1$, and $h(-1) = h(1) = 0$. Hence, for $\varepsilon > 0$, there are $\delta_1, \delta_2 \in (0, 1)$ such that $|h(x)| < \varepsilon$ holds for $x \in [-1, -\delta_1) \cup (\delta_2, 1]$. On the other hand, for $x \in [-\delta_1, \delta_2]$, we have

$$h(x) = f(x) - q_n^{[2]}(x) = g^{[2]}(z) - p_n^{[2]}(z) = \frac{\ell(z)}{(n+1)!} \left(g^{[2]}(z)\right)^{(n+1)}\Big|_{z=\xi}$$

for some $\xi \in (\lambda_1, \lambda_2)$ depending on $x$, with $\lambda_1 = \kappa_2(\delta_2)$ and $\lambda_2 = \kappa_2(-\delta_1) = -\kappa_2(\delta_1)$, in which

$$\begin{aligned} \ell(z) &= \prod_{i=0}^{n} \left(z - z_i^{[2]}\right) = \prod_{i=0}^{n} \left[z - \cos\left(\frac{i}{n}\pi\right)\right] \\ &= (z-1)(z+1) \prod_{i=1}^{n-1} \left[z - \cos\left(\frac{i}{n}\pi\right)\right] \\ &= \frac{z^2 - 1}{2^{n-1}} \cdot \frac{\sin(n \arccos z)}{\sin(\arccos z)} = -\frac{1}{2^{n-1}} \sin(n \arccos z) \sin(\arccos z), \end{aligned}$$

in view of $1 - z^2 = \sin^2(\arccos z)$ and the formulae with respect to the second-kind Chebyshev polynomial [13, pp. 4 and 156]. This combined with $|h(x)| < \varepsilon$ holding for $x \in [-1, -\delta_1) \cup (\delta_2, 1]$ gives (3.3) and completes the proof. □

By direct calculations, $q_n^{[1]}(x)$ and $q_n^{[2]}(x)$ given in Theorem 1 and Theorem 2 can be explicitly expressed as

$$\begin{aligned} q_n^{[1]}(x) &= p_n^{[1]}(\kappa_1(x)) = \sum_{k=0}^{n} \delta_k^{[1]} c_k^{[1]} T_k(\kappa_1(x)) = \sum_{k=0}^{n} \delta_k^{[1]} c_k^{[1]} \cos\left(k \arccos\left(\kappa_1(x)\right)\right) \\ &= \sum_{k=0}^{n} \delta_k^{[1]} c_k^{[1]} \cos\left(k \arccos\left(-\sin\left(\frac{n\pi}{2(n+1)} x\right)\right)\right) \\ &= \sum_{k=0}^{n} \delta_k^{[1]} c_k^{[1]} \cos\left(\frac{k\pi}{2} + \frac{kn\pi}{2(n+1)} x\right), \end{aligned} \tag{3.4}$$

$$\begin{aligned} q_n^{[2]}(x) &= p_n^{[2]}(\kappa_2(x)) = \sum_{k=0}^{n} \delta_k^{[2]} c_k^{[2]} T_k(\kappa_2(x)) = \sum_{k=0}^{n} \delta_k^{[2]} c_k^{[2]} \cos\left(k \arccos\left(\kappa_2(x)\right)\right) \\ &= \sum_{k=0}^{n} \delta_k^{[2]} c_k^{[2]} \cos\left(k \arccos\left(-\sin\left(\frac{\pi}{2} x\right)\right)\right) \\ &= \sum_{k=0}^{n} \delta_k^{[2]} c_k^{[2]} \cos\left(\frac{k\pi}{2} + \frac{k\pi}{2} x\right), \end{aligned} \tag{3.5}$$

respectively, with

$$c_k^{[1]} = \frac{2}{n+1} \sum_{i=0}^{n} y_i \cos\left(\frac{k(2i+1)}{2(n+1)}\pi\right), \quad \delta_k^{[1]} = \begin{cases} 1/2, & k = 0 \\ 1, & k = 1, \cdots, n \end{cases}$$

$$c_k^{[2]} = \frac{2}{n} \sum_{i=0}^{n} \delta_i^{[2]} y_i \cos\left(\frac{ki}{n}\pi\right), \qquad \delta_k^{[2]} = \begin{cases} 1/2, & k = 0, n \\ 1, & k = 1, \cdots, n-1 \end{cases}$$

On the other hand, in view of the uniqueness of polynomial interpolation, $q_n^{[1]}(x)$ and $q_n^{[2]}(x)$ can also be equivalently induced from (2.2) as follows:

$$q_n^{[1]}(x) = \left[\sum_{i=0}^{n} \frac{(-1)^{i+1} y_i \sin\left(\frac{2i+1}{2n+2}\pi\right)}{\sin\left(\frac{n\pi}{2n+2}x\right) + \cos\left(\frac{2i+1}{2n+2}\pi\right)}\right] \Bigg/ \left[\sum_{i=0}^{n} \frac{(-1)^{i+1} \sin\left(\frac{2i+1}{2n+2}\pi\right)}{\sin\left(\frac{n\pi}{2n+2}x\right) + \cos\left(\frac{2i+1}{2n+2}\pi\right)}\right],$$

$$q_n^{[2]}(x) = \frac{\dfrac{-y_0}{2\left[\sin\left(\frac{\pi}{2}x\right)+1\right]} + \displaystyle\sum_{i=1}^{n-1} \frac{(-1)^{i-1} y_i}{\sin\left(\frac{\pi}{2}x\right) + \cos\left(\frac{i}{n}\pi\right)} + \dfrac{(-1)^{n-1} y_n}{2\left[\sin\left(\frac{\pi}{2}x\right)-1\right]}}{\dfrac{-1}{2\left[\sin\left(\frac{\pi}{2}x\right)+1\right]} + \displaystyle\sum_{i=1}^{n-1} \frac{(-1)^{i-1}}{\sin\left(\frac{\pi}{2}x\right) + \cos\left(\frac{i}{n}\pi\right)} + \dfrac{(-1)^{n-1}}{2\left[\sin\left(\frac{\pi}{2}x\right)-1\right]}}.$$

### 3.3. *Method Naming*

Before naming the new method, we first observe what changes the two transformations defined by (3.1) and (3.2) have made to the function $f(x)$. For this purpose, we select eight functions from the literature that are widely used for measuring interpolation performance, and add two variants of one of these functions. These ten functions are listed in Table 1 and drawn in Figure 3 and Figure 4 (thick black solid lines). Among them, the regions with relatively rich variations are located either in the middle of $[-1, 1]$, near the two endpoints, or scattered throughout the entire interval. Therefore, the selection of these ten functions comprehensively represents the various difficulties in interpolation.

Table 1: Functions used in the numerical experiment.

| $i$ | $f_i(x)$ | Ref. | Description |
|---|---|---|---|
| 1 | $\dfrac{1}{1+25x^2}$ | [4] | The standard Runge function |
| 2 | $\tanh(10x)$ | [6] | Hyperbolic tangent function |
| 3 | $\dfrac{1}{1+1000(x+0.5)^2} + \dfrac{1}{\sqrt{1+1000(x-0.5)^2}}$ | [6] | Bimodal function |
| 4 | $\lvert x\rvert + 0.5x - x^2$ | [2] | Asymmetric hump: an attempt |
| 5 | $\mathrm{e}^{-5x^2}\sin(20x)$ | [16] | Gaussian-modulated oscillation |
| 6 | $\dfrac{\sin(50x)}{1+25x^2}$ | [6] | High-frequency oscillation |
| 7 | $\dfrac{\sin(30x)}{1+10x^2} \cdot \sin\left(\dfrac{10}{1+5x^2}\right)$ | [6] | Rational oscillation |
| 8 | $\cos^{19}(4\pi x)$ | [17] | Higher-order polynomial of cosine |
| 9 | $\dfrac{1}{1+1000(x+0.2)^2} + \dfrac{1}{\sqrt{1+1000(x-0.2)^2}}$ | [6] | A contracting variant of $f_3(x)$ |
| 10 | $\dfrac{1}{1+1000(x+0.8)^2} + \dfrac{1}{\sqrt{1+1000(x-0.8)^2}}$ | [6] | A dilating variant of $f_3(x)$ |

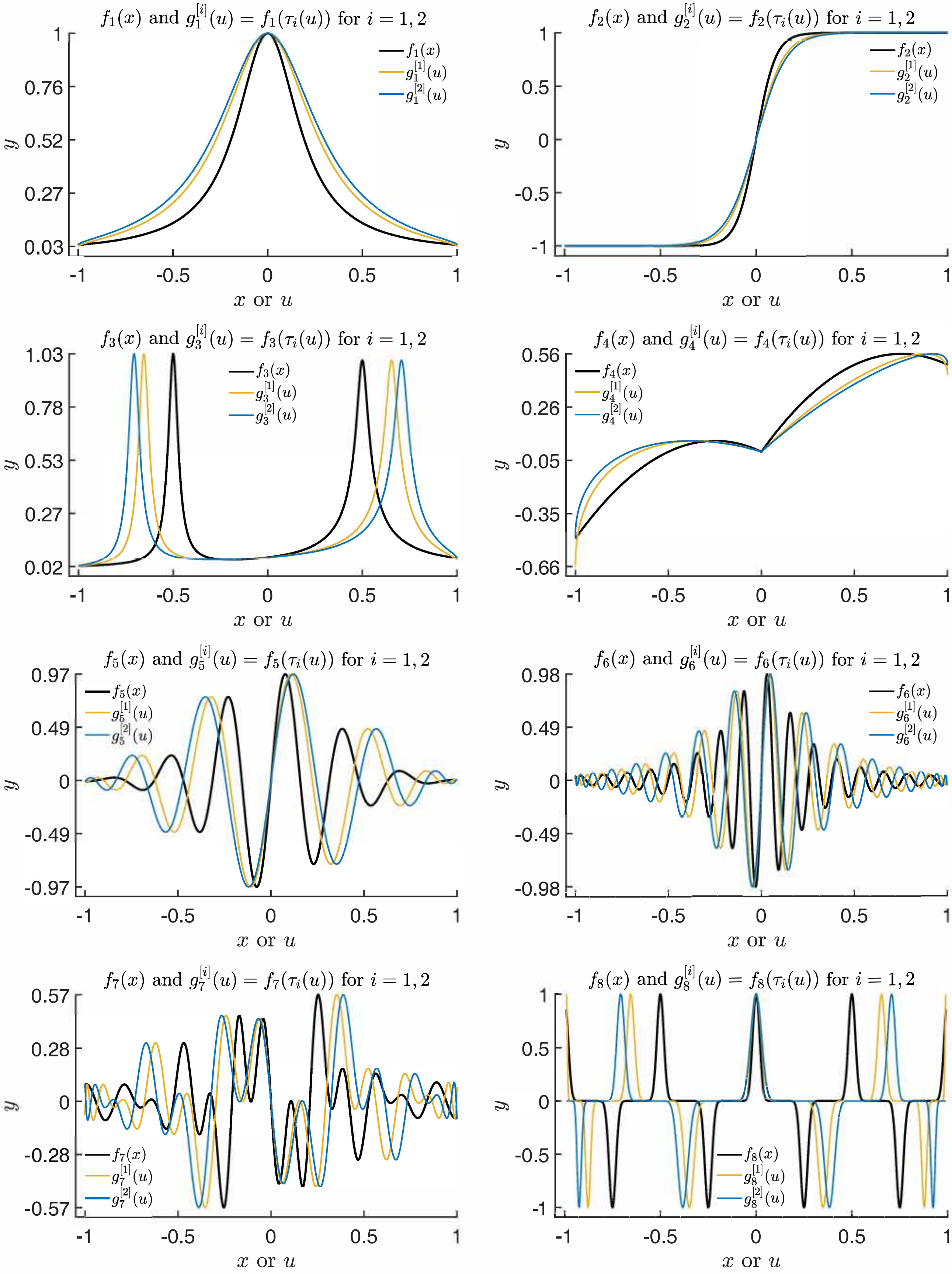

Figure 3: Functions selected from the literature used to measure interpolation performance.

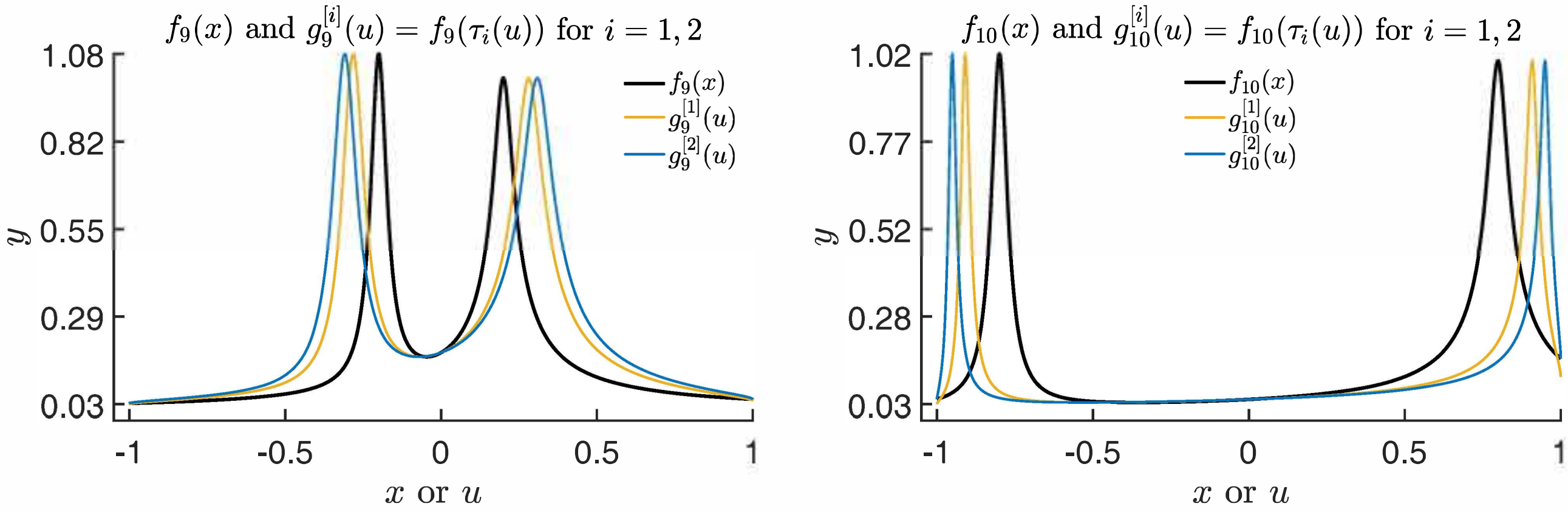

Figure 4: The two variants of the third function, one is contracting and the other is dilating.

For each of the ten functions $f_i(x)$, the transformed versions $g_i^{[1]}$ and $g_i^{[2]}$ are also shown in Figure 3 and Figure 4. The results indicate that both $g_i^{[1]}$ and $g_i^{[2]}$ can be interpreted as new configurations, shaped by a symmetric, wave-like force propagating from the center of $f_i(x)$ toward the endpoints. In short, this process resembles a wave-induced symmetric deformation. Consequently, we name this method "*symmetric wave interpolation* (SWI)".

Conceptually, SWI redistributes equidistant points by applying a symmetric central force, creating a non-uniformly compressed distribution for stable interpolation. After computations, the points revert to their initial equidistant layout, and the solution is mapped back, yielding high accuracy without the Runge phenomenon. Theoretically, SWI possesses the same computational complexity as that of Chebyshev interpolation, $O(n \ln n)$; Computationally, SWI may be slightly faster as it does not require the calculation of arccosine.

In order to distinguish between the two versions of SWI given in (3.4) and (3.5), we refer to them as *SWI of the first kind on equidistant points* (abbreviated *SWI1-Equid*) and *SWI of the second kind on equidistant points* (abbreviated *SWI2-Equid*), respectively. Accordingly, the *Chebyshev interpolation on Chebyshev points of the first kind* and *Chebyshev interpolation on Chebyshev points of the second kind*, which are denoted by (2.3) and (2.4), will also be abbreviated as *CI-Cheby1* and *CI-Cheby2* when necessary.

### *3.4. A preliminary observation*

To preliminarily observe the interpolation performance of SWI, we apply *SWI1-Equid* and *SWI2-Equid* to interpolate the Runge function $f_1(x)$ based on $(n + 1)$ equidistant points, where $n$ takes the same values as in Figure 2. The results are presented in Figure 5.

As shown in the figure, SWI effectively resists the Runge phenomenon while only using the conveniently available equidistant points, demonstrating its successful inheritance of the stability of Chebyshev interpolation. Additionally, we note from Figure 2 and Figure 5 that *SWI1-Equid* passes precisely through the endpoints, while *CI-Cheby1* not only fails to do so but also exhibits slight but non-negligible oscillations near the boundaries. This means that SWI has the potential to surpass Chebyshev interpolation in certain situations. In the next section, we will explore this possibility by means of all of the ten benchmark functions listed in Table 1.

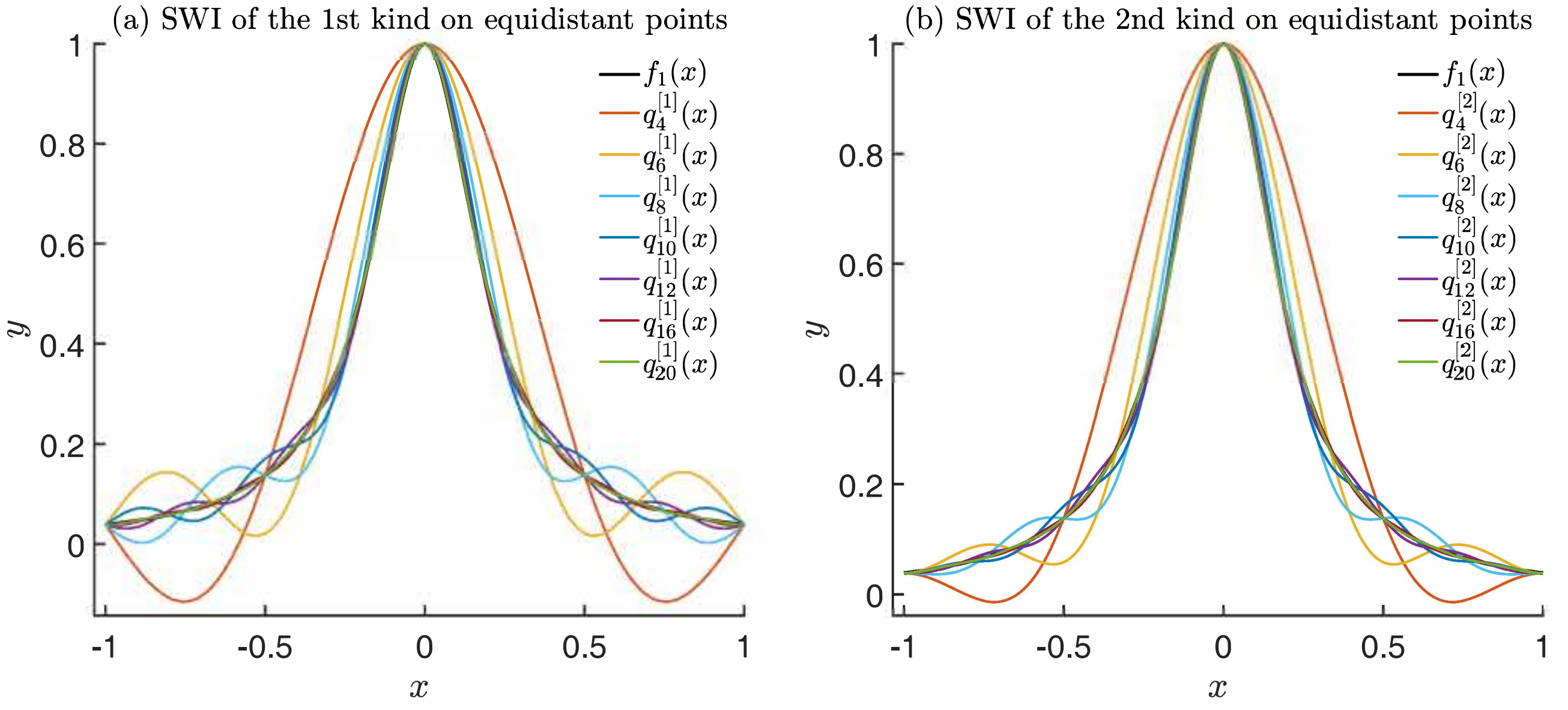

Figure 5: SWI of the first and second kinds of $f_1(x)$ based on the $(n+1)$ equidistant points for $n = 4, 6, 8, 10, 12, 16, 20$.

### *3.5. Robustness*

In practice, data collection often cannot achieve extremely high precision, as measurements are typically accurate to two or three decimal places, rather than a dozen or more. This introduces observational or round-off or other errors into the data. For convenience, here we consider the rounding errors of data.

To evaluate how such errors affect the robustness of SWI in real applications, we round the original interpolation data to only two significant digits, use SWI to interpolate $f_1(x)$ based on 13 equidistant points, and observe their impact numerically. Figure 6 presents the corresponding results. It demonstrates that even when the input data is rounded to only two significant digits, the results of SWI still closely match that obtained from the exact data. This observation demonstrates the robustness of SWI in the presence of data imprecision.

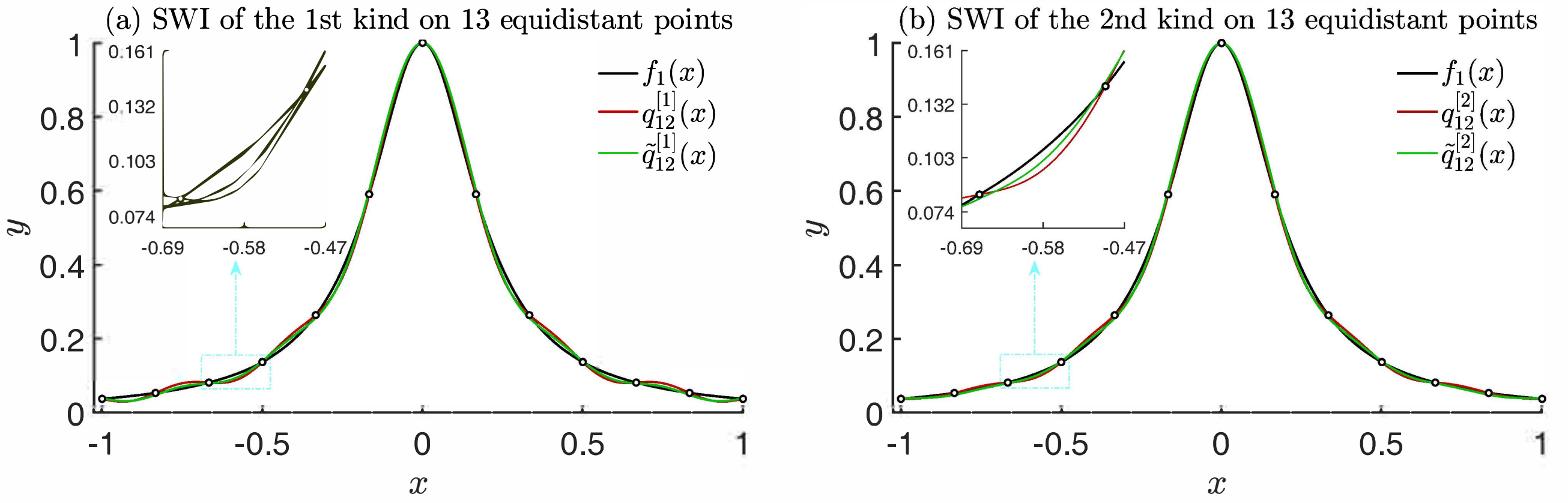

Figure 6: SWI of the first and second kinds of $f_1(x)$ based on the 13 equidistant points: $q_{12}^{[i]}(x)$ uses exact sample data, while $\tilde{q}_{12}^{[i]}(x)$ uses the data rounded to only two significant digits.

## 4. Numerical experiment

The preliminary test in Section 3 has shown that SWI successfully combines the *practical convenience of equidistant points* with the *numerical stability of Chebyshev interpolation on Chebyshev points*, indicating its potential to outperform Chebyshev interpolation in some cases. This section will conduct a comprehensive numerical experiment to systematically evaluate the advantages of the two SWI methods, while also analyzing how they work to identify the possible limitations to guide future improvements.

### *4.1. Experimental setup*

Before proceeding, we describe the following aspects of the experimental setup:

- *Interpolation methods*: In view of the numerical stability and the uniform convergence of Chebyshev interpolation [15], besides *SWI1-Equid/SWI2-Equid*, the experiment focuses solely on the two kinds of Chebyshev interpolation, *CI-Cheby1/CI-Cheby2*. The aim is to demonstrate SWI1-Equid/SWI2-Equid, while using only equidistant points, can achieve performance comparable to that of *Chebyshev points based* CI-Cheby1/CI-Cheby2.
- *Benchmark functions*: A total of eight benchmark functions widely used in the literature, plus two variants of the third function, are selected for the assessment of interpolation performance. The ten functions are presented in Table 1. As illustrated in Figure 3 and Figure 4, these functions cover a wide range of difficulties in interpolation.
- *Performance Metrics*: To select appropriate metrics for evaluating interpolation quality, we revisit the results of interpolating the Runge function $f_1(x)$ using *SWI1-Equid/SWI2-Equid* and *CI-Cheby1/CI-Cheby2*, as shown in Figure 2 and Figure 5. For a clearer observation of the underlying behavior, we extract the case of $n = 12$ and present it in Figure 7. By the figure, $q_{12}^{[j]}$ seems to be better than $p_{12}^{[j]}$ as it is closer to $f_1(x)$. Here, the "closeness" between $f_k(x)$ and the interpolant function, denoted by $\varphi(x)$, can be quantified using the following two types of errors:

  ◇ Maximal error: $\|f_k - \varphi\|_\infty \triangleq \max\limits_{x\in[-1,1]} |f_k(x) - \varphi(x)|$;

  ◇ Cumulative error: $\|f_k - \varphi\|_1 \triangleq \int_{-1}^{1} |f_k(x) - \varphi(x)|\, \mathrm{d}x$.

  Clearly, the first metric corresponds to the maximal vertical distance between $f_k(x)$ and $\varphi(x)$, while the second corresponds to the area between the curves bounded by $x = \pm 1$. In both cases, smaller error values indicate better interpolation performance.

### *4.2. Numerical results*

Based on the setup outlined above, the numerical experiment is carried out as follows. For each function in Table 1 and a given positive integer $n$, execute the following steps:

- *Chebyshev interpolation*: Perform CI-Cheby1/CI-Cheby2 on $(n + 1)$ Chebyshev points;
- *Equidistant interpolation*: Perform SWI1-Equid/SWI2-Equid on $(n+1)$ equidistant points;
- *Error calculation*: Compute both types of errors for each interpolation method.

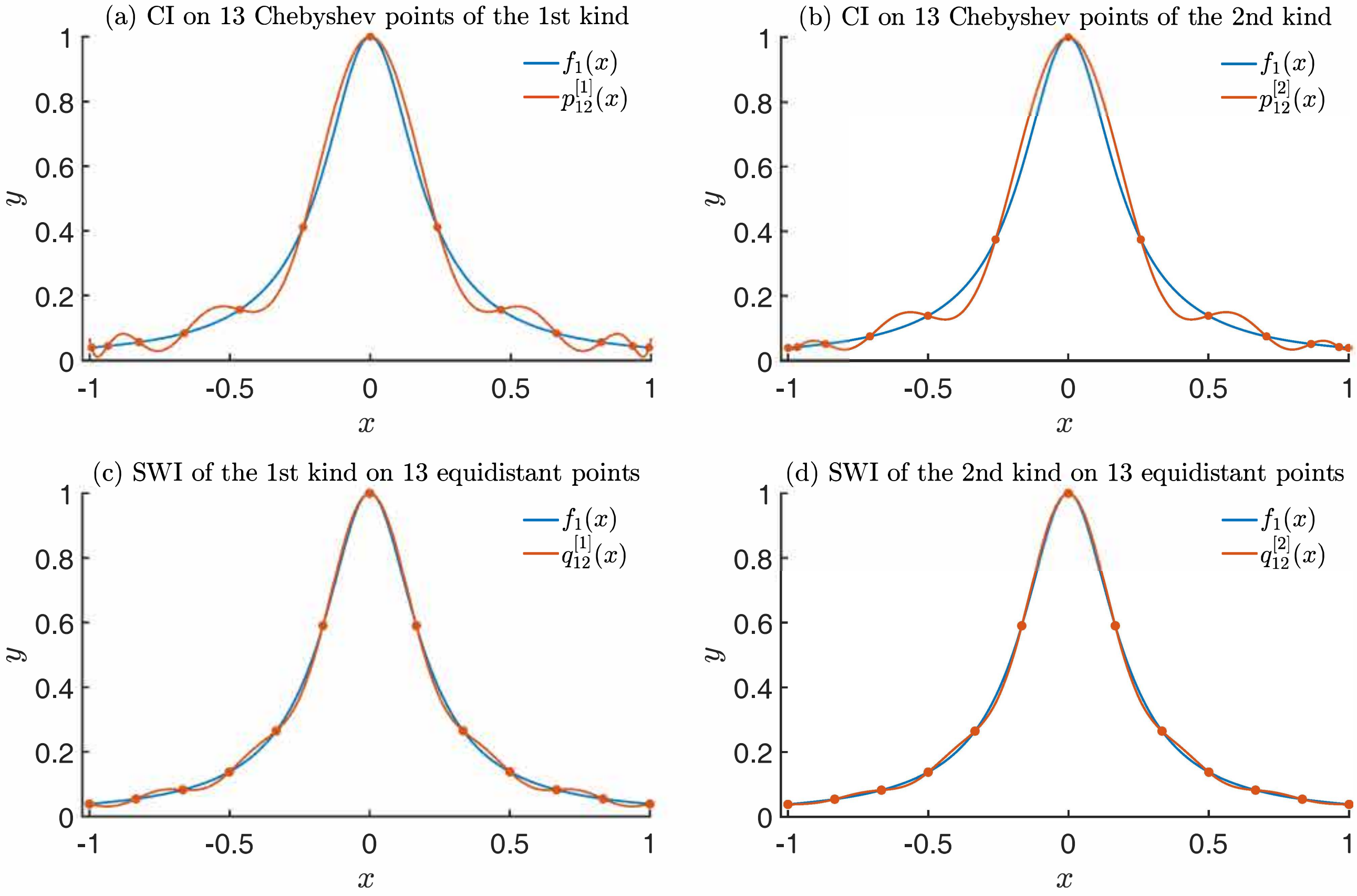

Figure 7: Chebyshev interpolation (CI) on 13 Chebyshev points and two kinds of SWI on 13 equidistant points for the Runge function $f_1(x)$.

Then, plot the computed errors against $n$. The results, displayed in Figure 8 and Figure A.1 through Figure A.9, lead to the following two major findings:

First, for all benchmark functions, SWI rapidly achieves a sufficiently small error level. This demonstrates that effective and stable global interpolation can indeed be accomplished using *readily available equidistant points*, without the need for *specially designed and hard-to-obtain Chebyshev points*. This represents a key contribution of this work.

Second, for almost all situations, SWI reaches a suitably small error level (e.g., not larger than 0.001) much faster than Chebyshev interpolation. In other words, SWI requires fewer equidistant points in most situations to match the interpolation quality that Chebyshev interpolation attains with more Chebyshev points.

To numerically show this conclusion, we list in Table 2 the minimal polynomial degree (i.e., the number of points minus one) required for Chebyshev interpolation and SWI across ten functions to first reach a *maximal error* or *cumulative error* below the sufficiently small positive scalar $\varepsilon$ (taken as 0.1 or 0.01 or 0.001), with that for each type taken as the smaller of the degrees required by its two variants (i.e., CI-Cheby1 and CI-Cheby2 for Chebyshev interpolation, and SWI1-Equid SWI2-Equid for SWI). This table provides more concrete evidences, clearly demonstrating the advantage of SWI over Chebyshev interpolation in terms of requiring fewer points to achieve the same quality of interpolation.

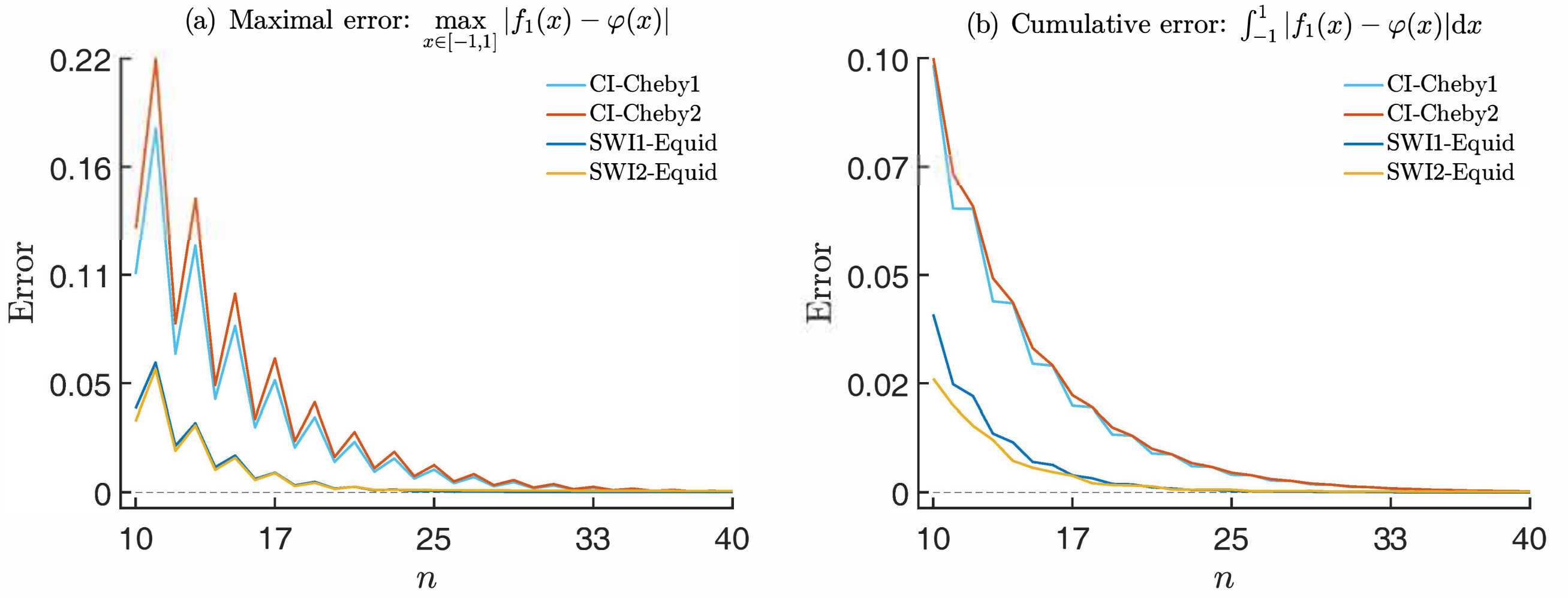

Figure 8: Maximal and cumulative errors versus $n$ (from 10 to 40) for interpolating $f_1(x)$ using *CI-Cheby1* ($\varphi = p_n^{[1]}$), *CI-Cheby2* ($\varphi = p_n^{[2]}$), *SWI1-Equid* ($\varphi = q_n^{[1]}$), and *SWI2-Equid* ($\varphi = q_n^{[2]}$).

Table 2: Minimal polynomial degree required for Chebyshev interpolation (CI) and SWI across ten functions to first reach a maximal error or cumulative error below $\varepsilon$.

| $i$ | Maximal error | | | | | | Cumulative error | | | | | |
|---|---|---|---|---|---|---|---|---|---|---|---|---|
| | $\varepsilon = 0.1$ | | $\varepsilon = 0.01$ | | $\varepsilon = 0.001$ | | $\varepsilon = 0.1$ | | $\varepsilon = 0.01$ | | $\varepsilon = 0.001$ | |
| | CI | SWI | CI | SWI | CI | SWI | CI | SWI | CI | SWI | CI | SWI |
| 1 | 12 | 8 | 24 | 16 | 34 | 24 | 9 | 6 | 21 | 14 | 33 | 22 |
| 2 | 17 | 11 | 31 | 21 | 45 | 29 | 12 | 8 | 28 | 18 | 42 | 27 |
| 3 | 66 | 48 | 129 | 96 | 192 | 140 | 20 | 18 | 80 | 58 | 140 | 106 |
| 4 | 6 | 4 | 60 | 40 | 596 | 380 | 5 | 4 | 15 | 13 | 49 | 39 |
| 5 | 25 | 15 | 29 | 17 | 33 | 19 | 25 | 15 | 29 | 17 | 33 | 19 |
| 6 | 59 | 37 | 69 | 45 | 81 | 195 | 56 | 36 | 68 | 44 | 80 | 53 |
| 7 | 45 | 31 | 59 | 37 | 65 | 122 | 45 | 29 | 54 | 37 | 65 | 43 |
| 8 | 98 | 56 | 140 | 88 | 182 | 120 | 83 | 56 | 131 | 88 | 165 | 120 |
| 9 | 78 | 50 | 147 | 95 | 218 | 140 | 22 | 14 | 93 | 61 | 157 | 104 |
| 10 | 44 | 50 | 89 | 99 | 133 | 148 | 13 | 14 | 56 | 63 | 96 | 105 |

### *4.3. Convergence of SWI*

Theorem 1 and Theorem 2 theoretically guarantee the convergence of *SWI1-Equid* and *SWI2-Equid* to the target function as the number of points, $n+1$, grows. To numerically validate this theoretical convergence, we examine two types of errors for SWI as $n$ increases from 100 to 1000. The results are shown in Figure 9 and Figure 10.

It is evident from the figures that both types of errors exhibit a decreasing trend and eventually approach zero as $n$ grows. These numerical results provide a support for our theoretical findings, collectively confirming the desirable convergence properties of SWI.

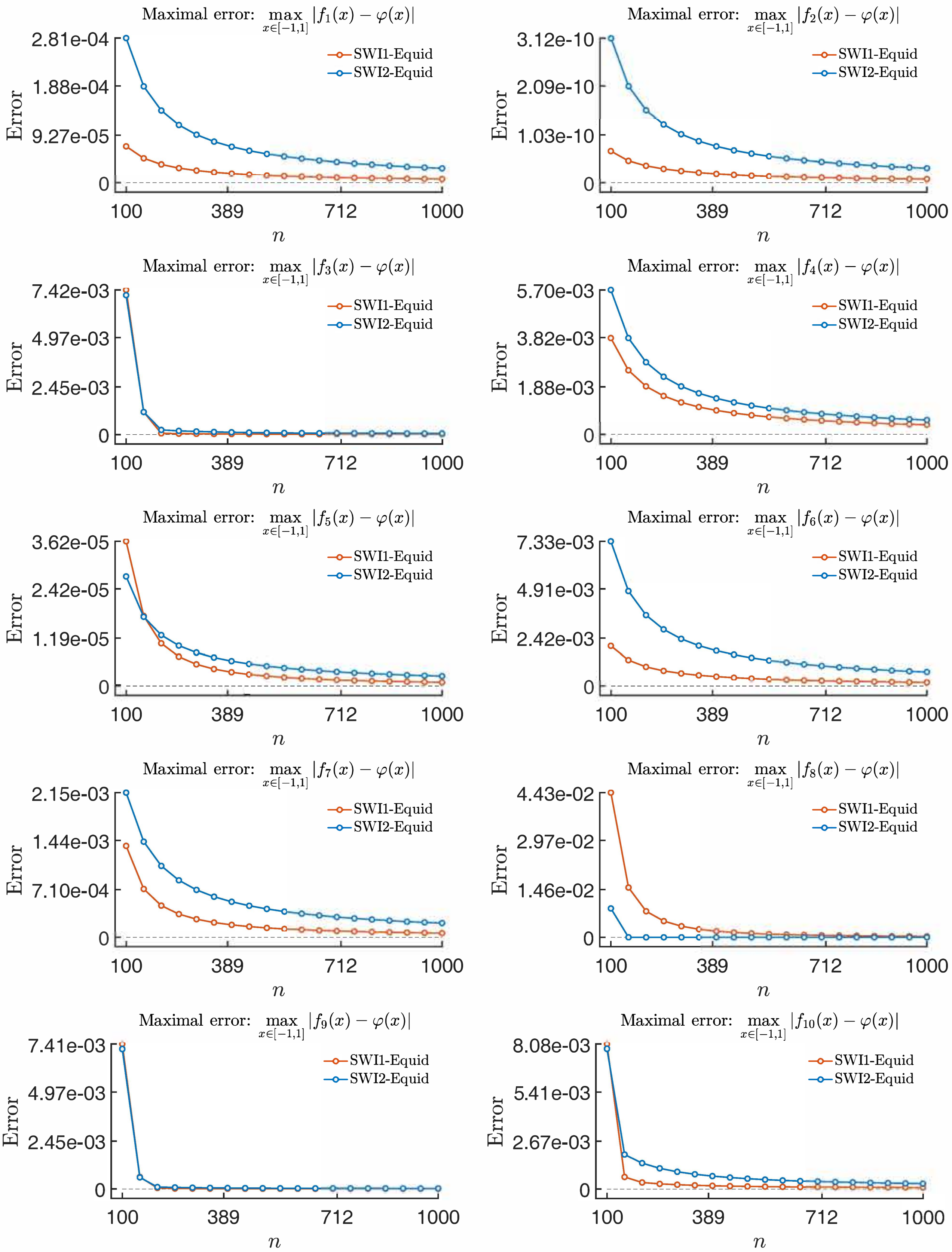

Figure 9: Maximal errors of SWI versus large $n$ (from 100 to 1000) with $\varphi = q_n^{[1]}$ or $q_n^{[2]}$.

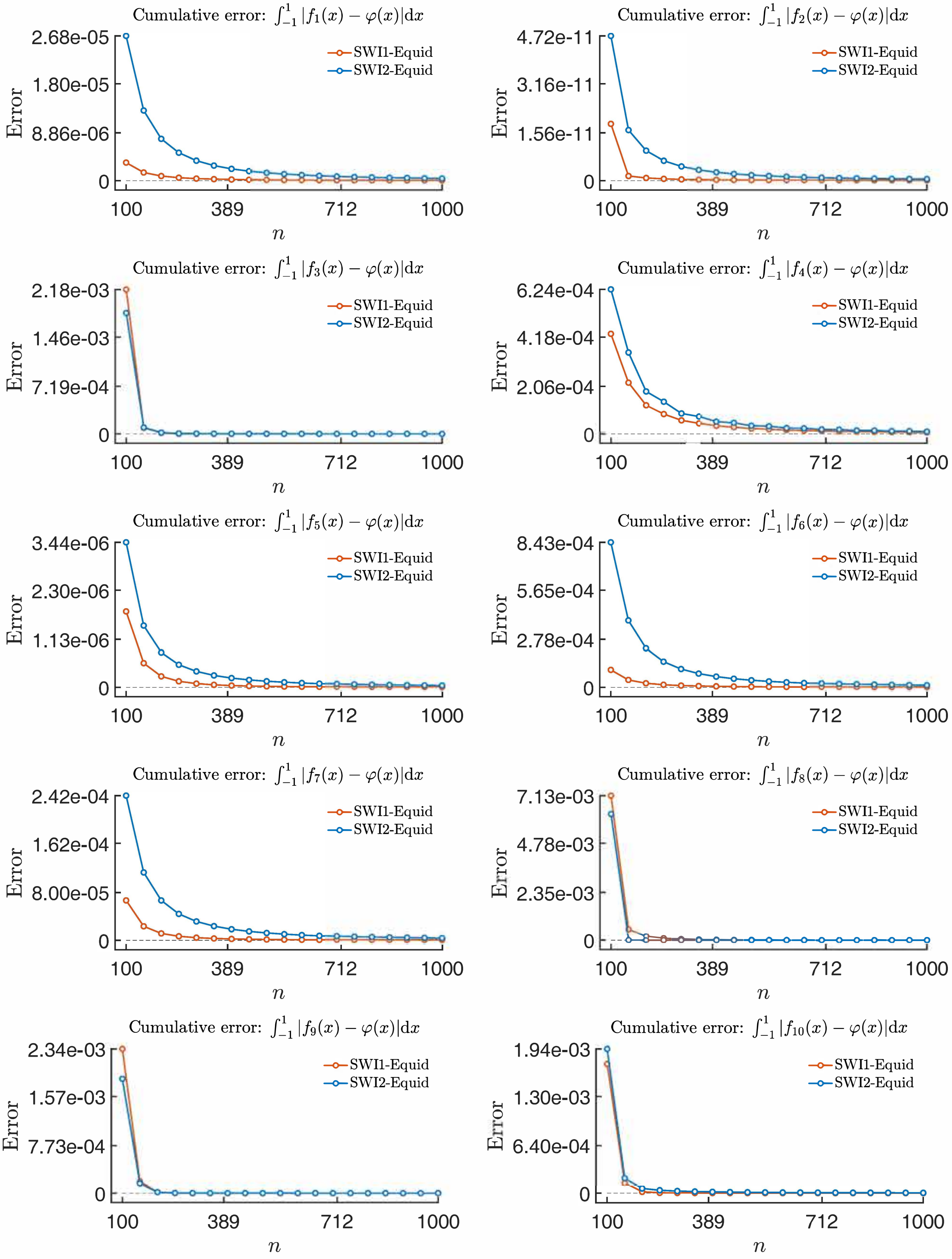

Figure 10: Cumulative errors of SWI versus large $n$ (from 100 to 1000) with $\varphi = q_n^{[1]}$ or $q_n^{[2]}$.

### 4.4. A limitation

Although the experiments demonstrate the excellent performance of SWI, it is important to point out a limitation. Specifically, both $\tau_1$ and $\tau_2$ exacerbate the difficulty of interpolating the transformed functions $g^{[1]}$ and $g^{[2]}$ at the endpoints compared to the original function $f$.

Table 2 also reveals a potential limitation of SWI, which we will analyze in depth to identify the root causes, assess its impact on the practical applications of SWI, and suggest a direction for future research. The potential limitation is that, for the tenth function $f_{10}(x)$, SWI consistently requires more equidistant points than Chebyshev interpolation to achieve comparable interpolation accuracy (in terms of both types of errors). To investigate this behavior, we examine the plots of $f_9(x)$ and $f_{10}(x)$ in Figure 4. Note that these two functions were intentionally derived from $f_3(x)$ to highlight this particular limitation of SWI.

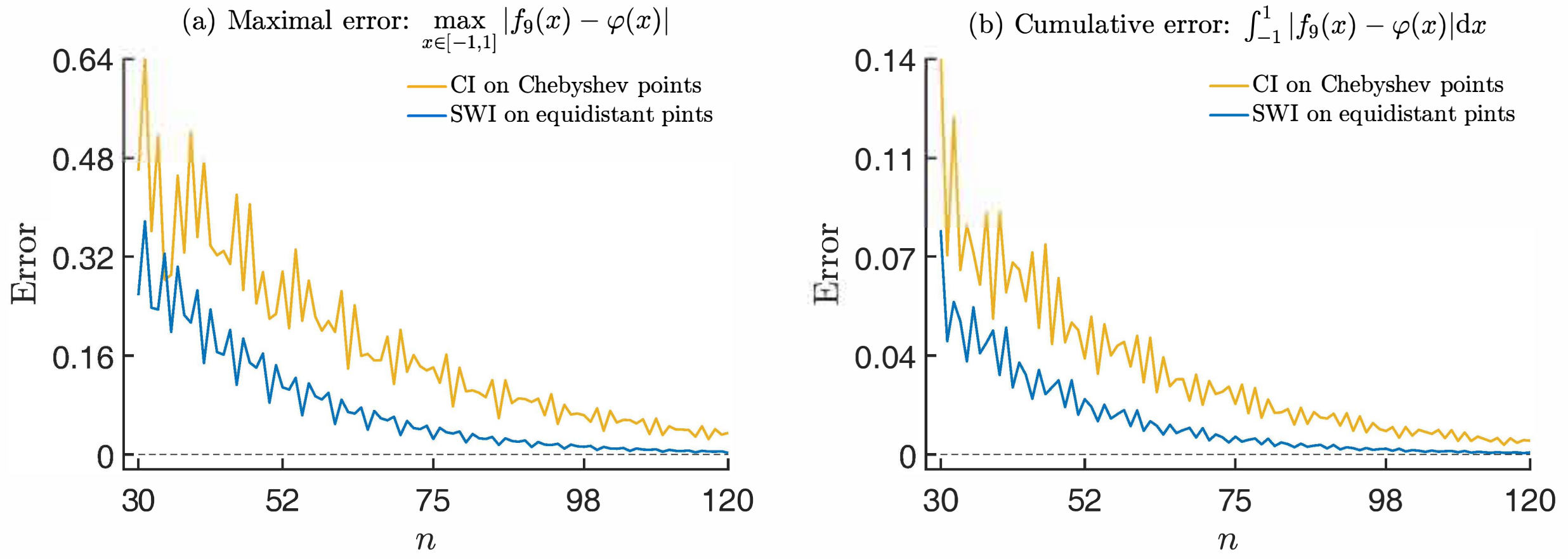

Figure 11: Maximal and cumulative errors of Chebyshev interpolation and SWI versus $n$ (from 30 to 120) for interpolating $f_9(x)$, in which $\varphi = \left(p_n^{[1]} + p_n^{[2]}\right)/2$ or $\left(q_n^{[1]} + q_n^{[2]}\right)/2$.

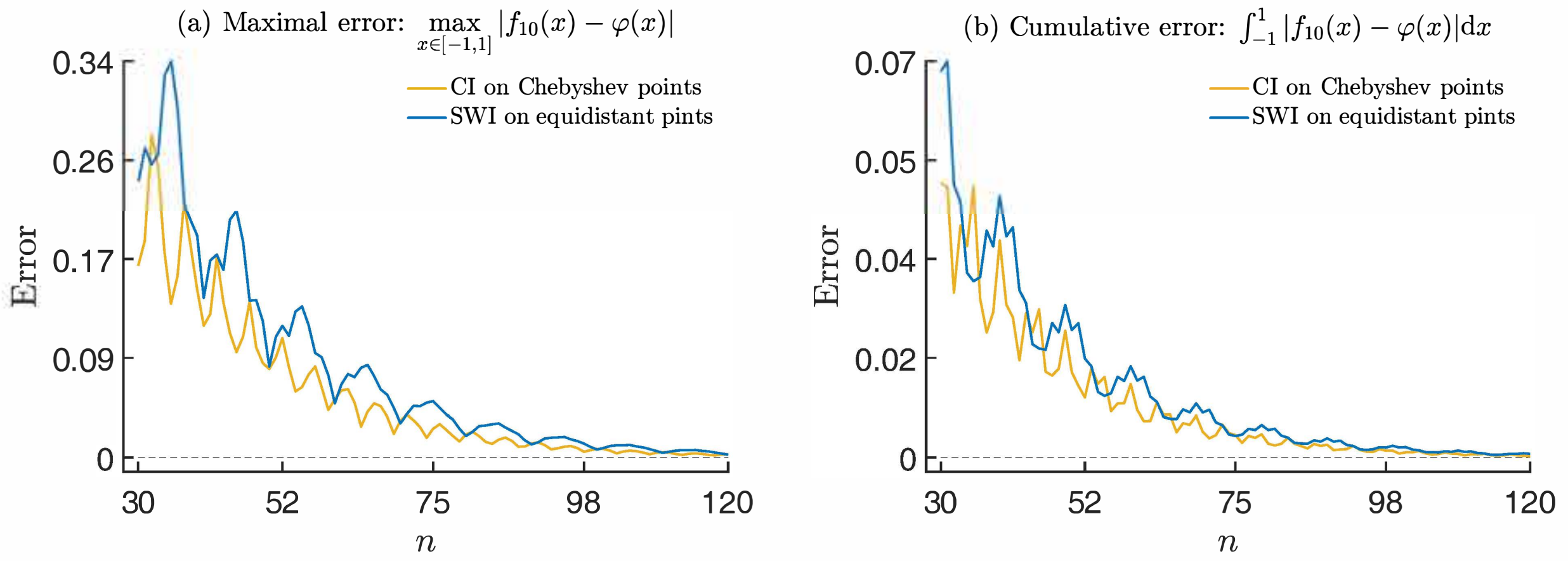

Figure 12: Maximal and cumulative errors of Chebyshev interpolation and SWI versus $n$ (from 30 to 120) for interpolating $f_{10}(x)$, in which $\varphi = \left(p_n^{[1]} + p_n^{[2]}\right)/2$ or $\left(q_n^{[1]} + q_n^{[2]}\right)/2$.

As shown in Figure 4, the left subplot of the ninth function $f_9(x)$ shows that the region of high variation is concentrated mainly in the central part of the interval $[-1, 1]$. After making the *symmetric wave mapping* $x = \tau_i(z)$, the high-varying part of $g_9^{[i]}(z)$ shifts away from the center and moves somewhat closer to the endpoints (though still remaining relatively distant) at which SWI performs well. In contrast, the right subplot of $f_{10}(x)$ shows that the high-varying region is already very close to the endpoints. After the same transformation, $x = \tau_i(z)$, this region moves even closer to both endpoints. Although SWI still performs well in this case, the interpolation task becomes noticeably more challenging, requiring more points to maintain high accuracy. Figure 11 and Figure 12 further illustrate how the two types of errors for Chebyshev interpolation and SWI evolve against $n$, confirming this observation.

This indicates that the shape of the function may influence the difficulty of interpolation using SWI. Specifically, when the region with high variation is closer to the center of the interval, SWI performs better; when it is closer to the endpoints, the interpolation accuracy tends to be affected. Given that SWI inherits traits from Chebyshev interpolation, we suspect this limitation may also exist in Chebyshev interpolation, though possibly amplified in SWI. To check this guess and measure the extent of the effect, we split the cumulative error into two parts: one from the endpoint regions (i.e., $[-1, -0.5]$ and $[0.5, 1]$) and the other from the central region (i.e., $[-0.5, 0.5]$). The results are shown in Figure 13 and Figure 14.

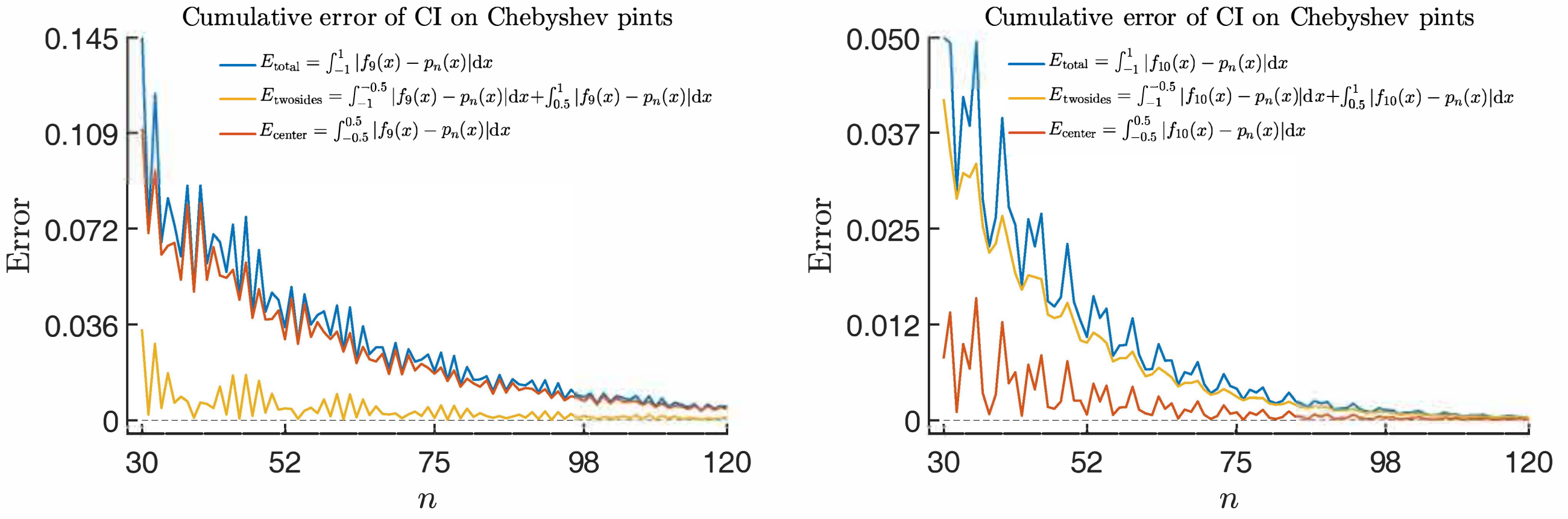

Figure 13: Partitioned cumulative errors of Chebyshev interpolation on Chebyshev points versus $n$ (from 30 to 120) for interpolating $f_9(x)$ and $f_{10}(x)$.

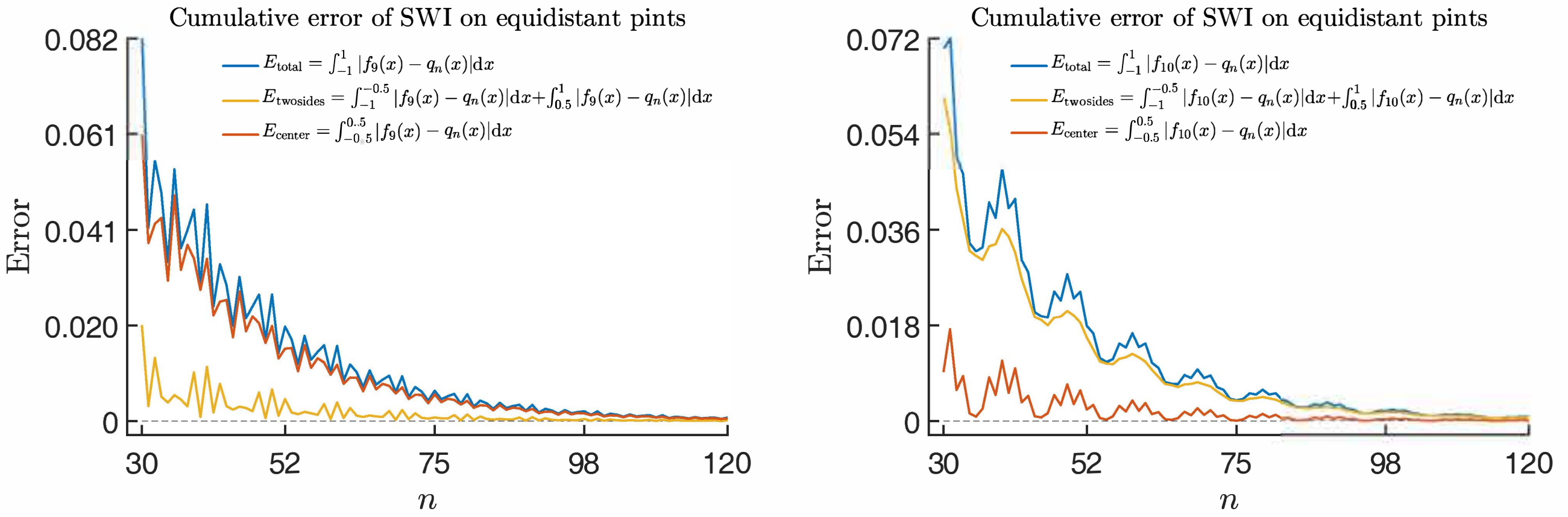

Figure 14: Partitioned cumulative errors of SWI on equidistant points versus $n$ (from 30 to 120) for interpolating $f_9(x)$ and $f_{10}(x)$.

The *left subplots* of these two figures indicate that, for both SWI and Chebyshev interpolation, the cumulative error is largely due to the central part of the interval $[-1, 1]$ when the function exhibits major variations near its center. Conversely, the *right subplots* of these two figures show that the error is dominated by the endpoint regions when the main variations occur there. Therefore, this limitation of SWI stems from two factors: one is the inherent property inherited from Chebyshev interpolation, and the other is the increased interpolation difficulty resulting from shifting points toward the endpoints. As Platte et al [5] showed, an interpolation method that simultaneously offers rapid convergence, high stability, and broad applicability remains fundamentally out of reach if we use only equidistant points. We can, however, make some efforts on this goal. Therefore, exploring new approaches to overcome this drawback will be a focus of our future work.

## 5. Concluding remarks

In this paper, we put forward the *symmetric wave interpolation* (SWI), a straightforward yet powerful technique for high-precision global interpolation on equidistant points. Our central finding shows that SWI effectively transfers the numerical stability of Chebyshev interpolation to the practical context of equidistant points. This achievement bridges a long-standing divide, successfully combining the practical convenience of equidistant points with the computational robustness previously exclusive only to specialized points like Chebyshev points. In essence, the method delivers a best-of-both-worlds outcome.

Experimental analyses reveal two key characteristics of SWI. First, its performance is relatively sensitive to the spatial distribution of a function's features (inheriting from Chebyshev interpolation): it excels when high-variation regions are centrally located but requires more points when these regions are near the boundaries. Second, although achieving very high precision can be less computationally efficient for SWI compared to Chebyshev interpolation, it is capable of reaching any specified accuracy level given adequate points. Theorem 1 and Theorem 2 offer a theoretical support for this finding. This represents a trade-off in efficiency rather than a fundamental limit on achievable accuracy.

In summary, this work makes a contribution to addressing the skepticism surrounding stable global interpolation that uses equidistant points. By demonstrating performance comparable to Chebyshev interpolation, SWI shows that global interpolation can be both robust and simple. The limitations of SWI (as well as that of Chebyshev interpolation), particularly the sensitivity to feature location, also provide a clear path for future work. We will focus on developing adaptive strategies or better point selection to make SWI more efficient and robust.

## Appendix A. Two types of errors for interpolating $f_2(x)$ through $f_{10}(x)$

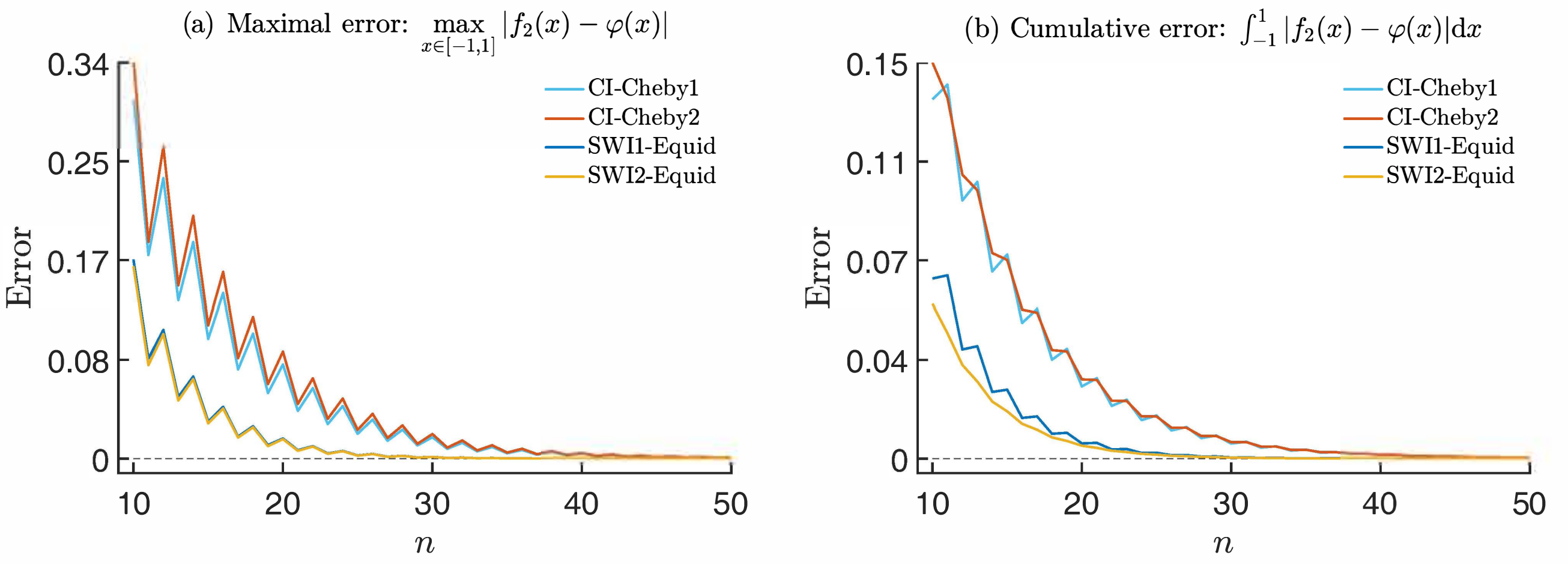

Figure A.1: Maximal and cumulative errors versus $n$ (from 10 to 50) for interpolating $f_2(x)$.

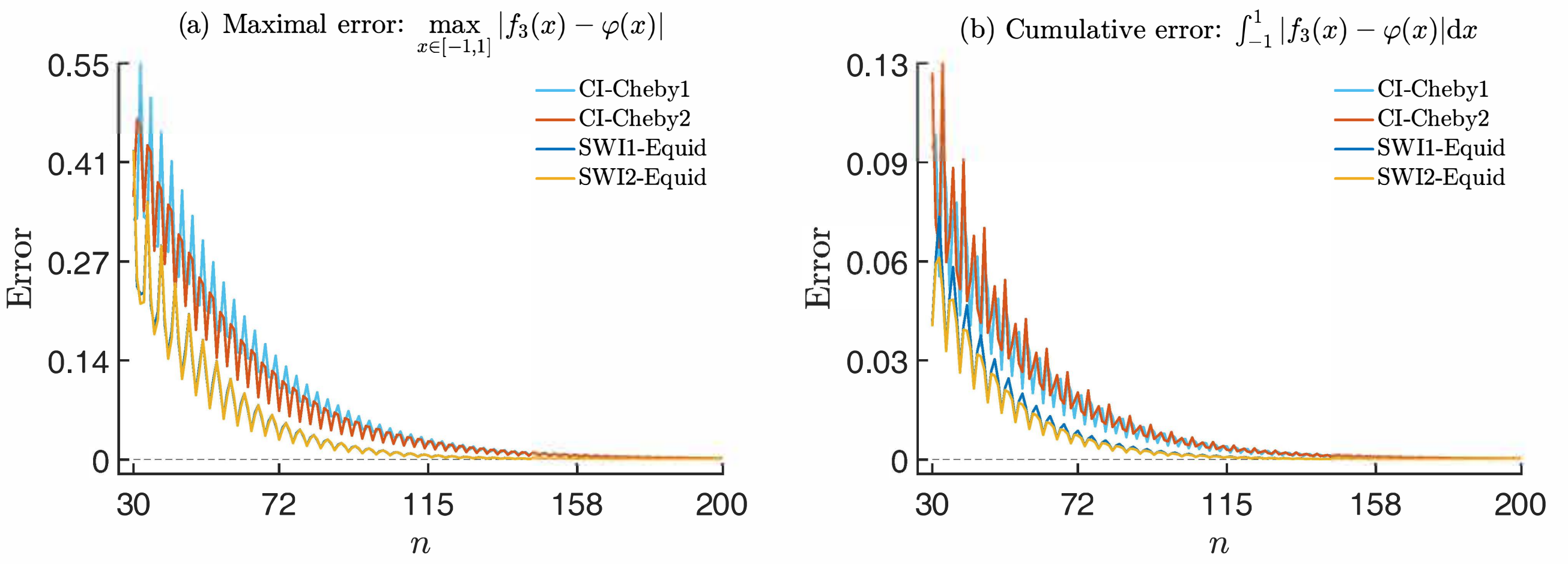

Figure A.2: Maximal and cumulative errors versus $n$ (from 30 to 200) for interpolating $f_3(x)$.

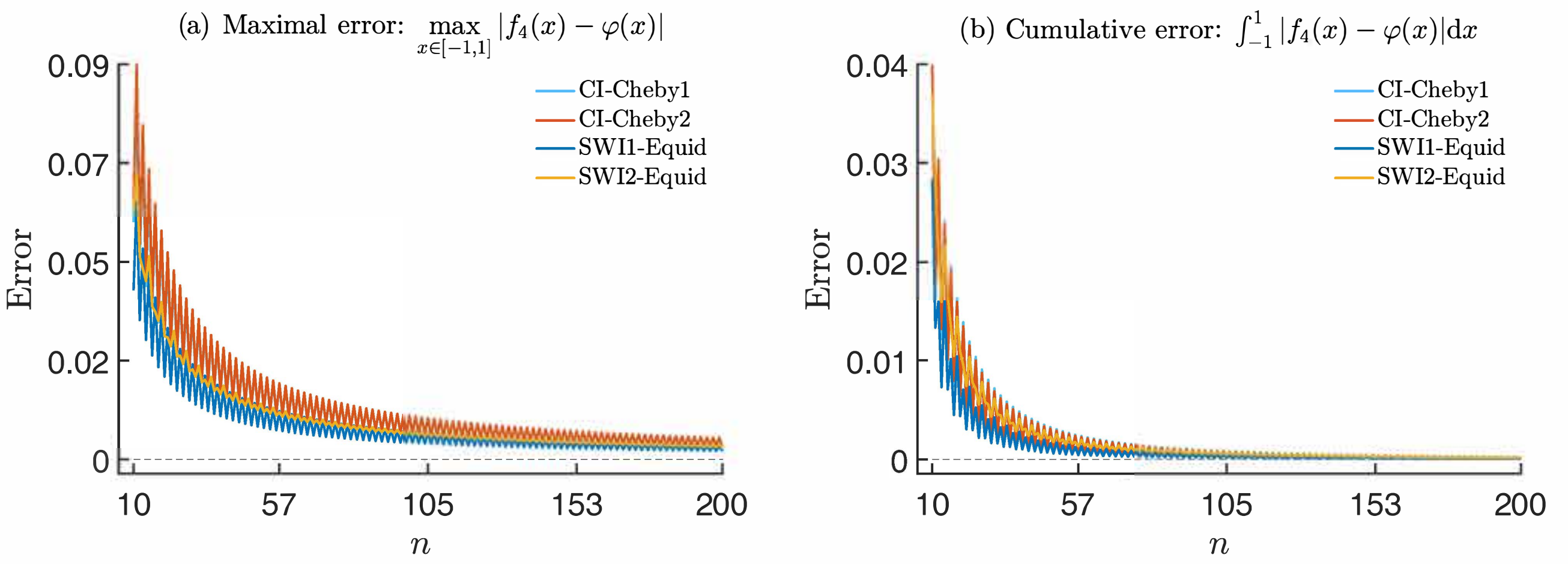

Figure A.3: Maximal and cumulative errors versus $n$ (from 10 to 200) for interpolating $f_4(x)$.

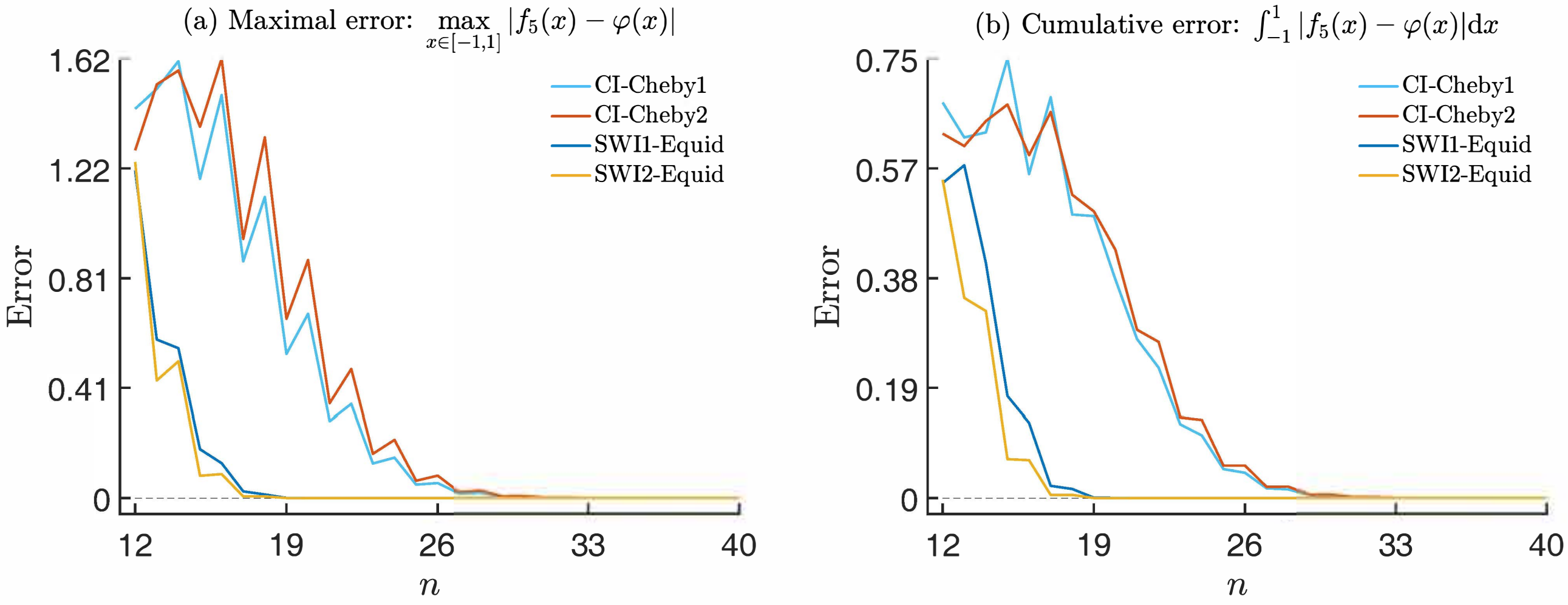

Figure A.4: Maximal and cumulative errors versus $n$ (from 12 to 40) for interpolating $f_5(x)$.

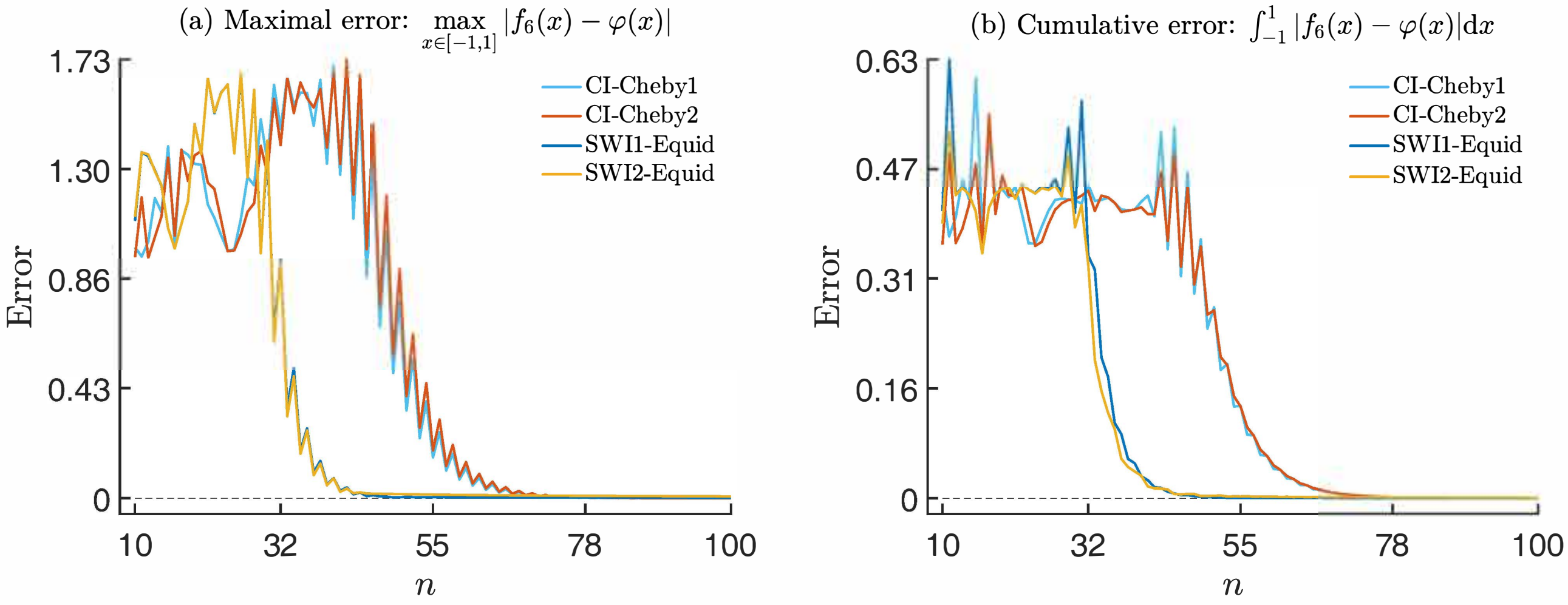

Figure A.5: Maximal and cumulative errors versus $n$ (from 10 to 100) for interpolating $f_6(x)$.

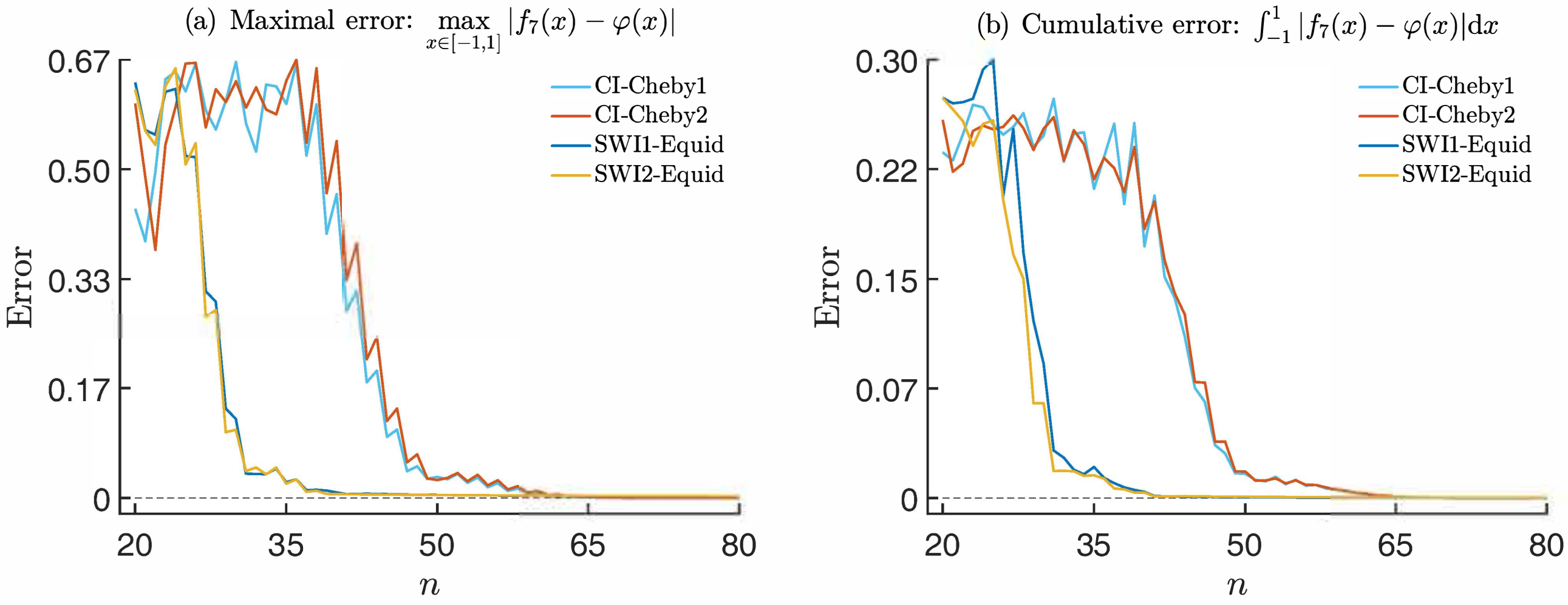

Figure A.6: Maximal and cumulative errors versus $n$ (from 20 to 80) for interpolating $f_7(x)$.

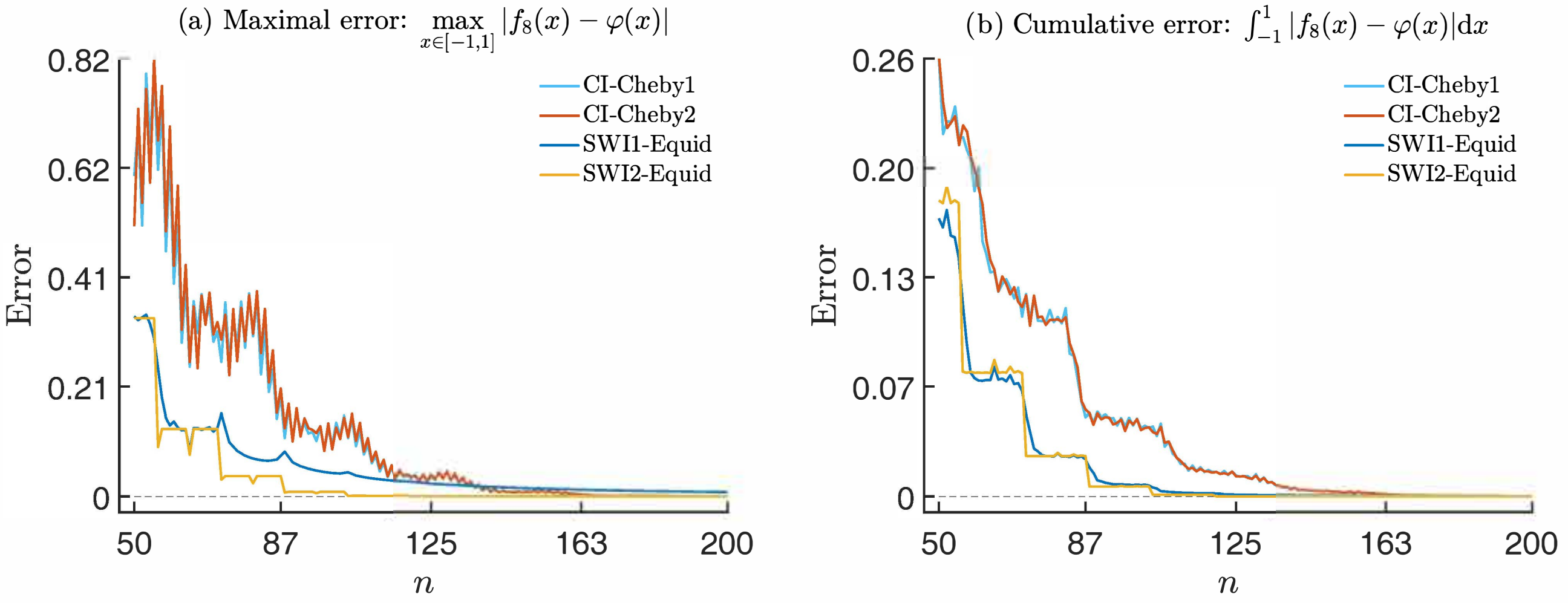

Figure A.7: Maximal and cumulative errors versus $n$ (from 50 to 200) for interpolating $f_8(x)$.

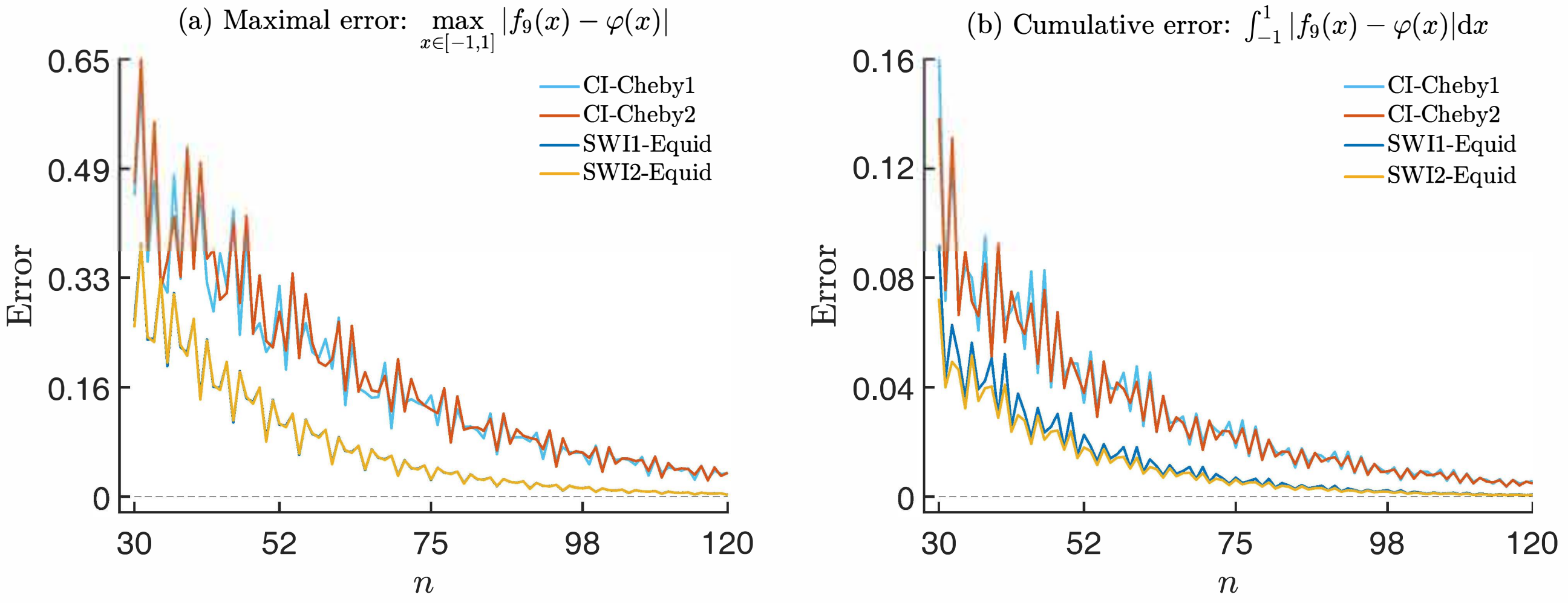

Figure A.8: Maximal and cumulative errors versus $n$ (from 30 to 120) for interpolating $f_9(x)$.

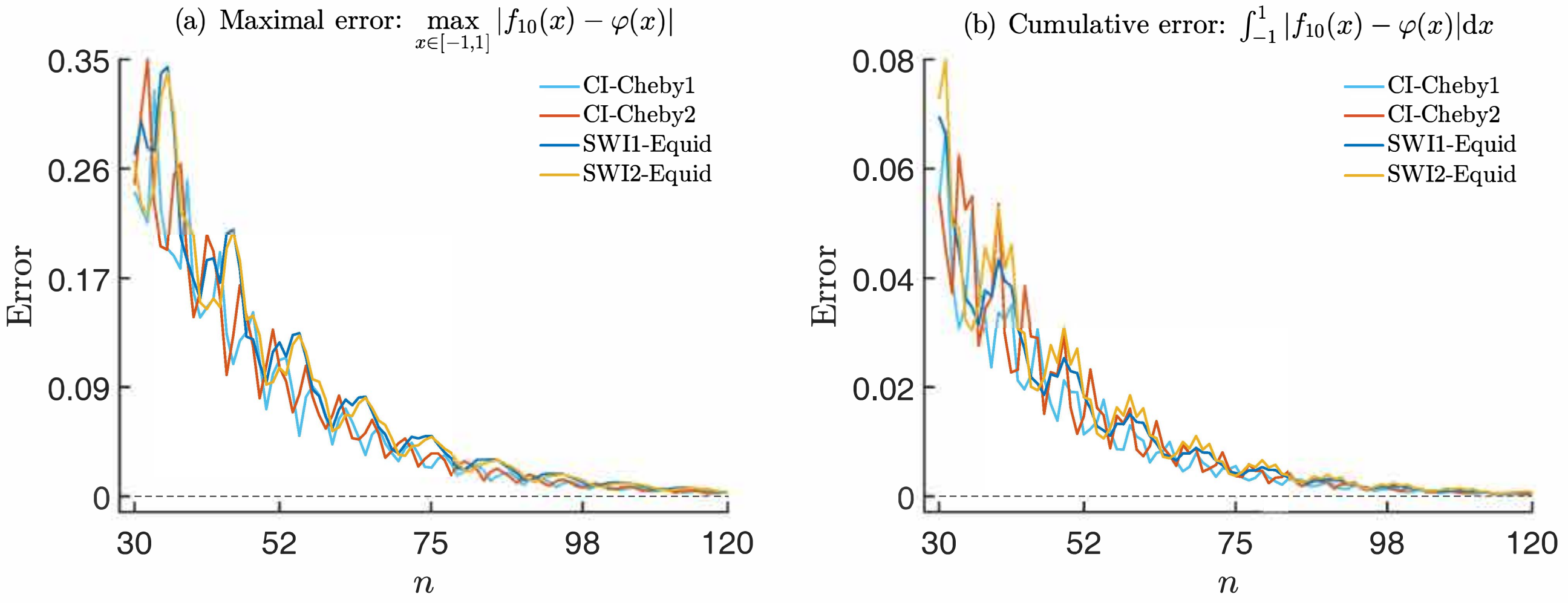

Figure A.9: Maximal and cumulative errors versus $n$ (from 30 to 120) for interpolating $f_{10}(x)$.